\documentclass[opre,nonblindrev]{informs3}

\OneAndAHalfSpacedXI % current default line spacing
\usepackage{endnotes}
\let\footnote=\endnote

\usepackage{dsfont}
\usepackage{algorithm}
\usepackage{algpseudocode}
\usepackage{subcaption}
\usepackage{multirow}
\usepackage[section]{placeins}
\newfloat{algorithm}{t}{lop}

\usepackage{natbib}
 \bibpunct[, ]{(}{)}{,}{a}{}{,}%
 \def\bibfont{\small}%
\TheoremsNumberedThrough     % Preferred (Theorem 1, Lemma 1, Theorem 2)
\ECRepeatTheorems

\EquationsNumberedThrough    % Default: (1), (2), ...
\begin{document}
%%%%%%%%%%%%%%%%

%\setcounter{page}{1}
%\DOI{10.1287/xxxx.0000.0000}%

\RUNAUTHOR{Mintz et al.}

\RUNTITLE{Behavioral Analytics for Myopic Agents}

% \TITLE{Bundling Information Goods of Decreasing Value}
\TITLE{Behavioral Analytics for Myopic Agents}

\ARTICLEAUTHORS{%
\AUTHOR{Yonatan Mintz, Anil Aswani, Philip Kaminsky}
\AFF{Department of Industrial Engineering and Operations Research, University of California, Berkeley, CA 94720, \EMAIL{\{aaswani,kaminsky,ymintz\}@berkeley.edu}}

\AUTHOR{Elena Flowers}
\AFF{Department of Physiological Nursing, School of Nursing, University of California, San Francisco, CA 94143, \EMAIL{	elena.flowers@ucsf.edu}}

\AUTHOR{Yoshimi Fukuoka}
\AFF{Department of Physiological Nursing \& Institute for Health \& Aging, School of Nursing, University of California, San Francisco, CA 94143, \EMAIL{	Yoshimi.Fukuoka@ucsf.edu}}

%\AUTHOR{Philip Kaminsky}
%\AFF{Department of Industrial Engineering and Operations Research, University of California, Berkeley, CA 94720, \EMAIL{kaminskyieor.berkeley.edu}}
%\AUTHOR{Yonatan Mintz}
%\AFF{Department of Industrial Engineering and Operations Research, University of California, Berkeley, CA 94720, \EMAIL{yminz@berkeley.edu}} %, \URL{}}
} % end of the block
\ABSTRACT{Many multi-agent systems have the structure of a single coordinator providing behavioral or financial incentives to a large number of agents. In such settings, two challenges faced by the coordinator are a finite budget from which to allocate incentives, and an initial lack of knowledge about the utility function of the agents. Here, we present a behavioral analytics approach to solve the coordinator's problem when the agents make decisions by maximizing utility functions that depend on prior system states, inputs, and other parameters that are initially unknown and subject to partially unknown temporal dynamics.  Our behavioral analytics framework involves three steps: first, we develop a behavioral model that describes the decision-making process of an agent; second, we use data to estimate behavioral model parameters for each agent and then use these estimates to predict future decisions of each agent; and third, we use the estimated behavioral model parameters to optimize a set of costly incentives to provide to each agent.  In this paper, we describe a specific set of tools, models, and approaches that fit into this framework, and that adapt models and incentives as new information is collected by repeating the second and third steps of this framework.  Furthermore, we prove that the incentives computed by this adaptive approach are asymptotically optimal with respect to a given loss function that describes the coordinator's objective.  We optimize incentives utilizing a decomposition scheme, where each sub-problem solves the coordinator's problem for a single agent, and the master problem is a pure integer program.  We conclude with a simulation study to evaluate the effectiveness of our behavioral analytics approach in designing personalized treatment plans for a weight loss program. The results show that our approach maintains efficacy of the program while reducing its costs by up to 60\%, while adaptive heuristics provide substantially less savings.

}%

%\KEYWORDS{deterministic inventory theory; infinite linear programming duality;
%  existence of optimal policies; semi-Markov decision process; cyclic schedule}

 %\HISTORY{This paper was
%first submitted on April 12, 1922 and has been with the authors for
%83 years for 65 revisions.}
\maketitle
%%%%%%%%%%%%%%%%%%%%%%%%%%%%%%%%%%%%%%%%%%%%%%%%%%%%%%%%%%%%%%%%%%%%%%

% Samples of sectioning (and labeling) in OPRE
% NOTE: (1) \section and \subsection do NOT end with a period
%       (2) \subsubsection and lower need end punctuation
%       (3) capitalization is as shown (title style).
%
%\section{Introduction.}\label{intro} %%1.
%\subsection{Duality and the Classical EOQ Problem.}\label{class-EOQ} %% 1.1.
%\subsection{Outline.}\label{outline1} %% 1.2.
%\subsubsection{Cyclic Schedules for the General Deterministic SMDP.}
%  \label{cyclic-schedules} %% 1.2.1
%\section{Problem Description.}\label{problemdescription} %% 2.
% Text of your paper here
\section{Introduction}

The increasing availability of data presents an opportunity to transform the design of incentives (i.e., costly inputs that are provided to agents to modify their behavior and decisions) from a single analysis into an adaptive and dynamic process whereby the incentive design is optimized as new data becomes available.  Historically, this adaptive setting has been studied under the framework of repeated games \citep{radner1985,fudenberg1994,laffont2002}, where researchers have focused on the analysis and identification of structural properties of effective policies, and on equilibria.   In contrast, continuing advances in optimization software and statistical estimation tools, utilized with the vast amount of data now available in many settings, enable a new approach that in many circumstances has the potential to lead to practical tools for designing effective incentives in real-world settings.   This approach, which we call \emph{behavioral analytics}, is built around a three step framework:  first, we develop a behavioral model that describes the decision-making process of an agent; next, we iterate repeatedly over two steps as new information is collected.  In the second step, we use data to estimate behavioral model parameters for each agent and then use these estimates to predict future decisions of each agent; and in the third, we use the estimated behavioral model parameters to optimize a set of costly incentives to provide to each agent.  In this paper, we describe a specific set of tools, models, and approaches that fit into this framework, and that adapt models and incentives as new information is collected while the second and third steps of the framework are repeated.  

% * <kaminsky@berkeley.edu> 2017-02-03T06:00:09.302Z:
%
% ^.

Specifically, we consider the following discrete-time setting: There is a large pool of agents each with a set of utility function parameters (which we will refer to as the motivational state) and system state at time $t$, and each agent makes a decision at $t$ by maximizing a myopic utility function.  A single coordinator makes noisy observations at $t$ of the system states and decisions of each agent, and then assigns behavioral or financial incentives (e.g., bonuses, payments, behavioral goals, counseling sessions) at $t$ to a subset of agents. The incentives change the motivational and system states of assigned agents at time $t+1$, while the motivational and system states of non-assigned agents evolve at $t+1$ according to some dynamics.  This process repeats, and time $t$ advances towards infinity in unit increments.  Here, the coordinator's problem is to decide what incentives to provide to which agents in order to minimize the coordinator's loss function, a function that depends on the system states and decisions of all agents.  This problem is challenging because the motivational states of agents are neither known nor measured by the coordinator, because agents make decisions by maximizing an unknown utility function, because measurements are noisy, and because the coordinator has a fixed budget (over a specified time horizon) from which to allocate incentives.

\subsection{Potential Applications for Behavioral Analytics}

The setting described above is found in many domains, including personalized healthcare, demand response programs, and franchise logistics.  Below, we elaborate on these potential applications of our framework.  The first application is the design of a weight loss program.  The coordinator is a clinician and the agents are individuals trying to lose body weight.  The next application is the design of a demand response program in which the coordinator is an electric utility company and the agents are homeowners (since they consume electricity).  The final application is in franchise logistics, where a parent company is a coordinator for a group of franchisees (who are the agents).

\subsubsection{Weight Loss Programs}
\label{sec:health_app}
In a clinically-supervised weight loss program, a clinician provides two types of behavioral incentives to a group of individuals who are trying to lose body weight.  The first type of incentive is behavioral goals provided to each individual by the clinician, and it is costless when communication costs are negligible, as is the case with mobile phone-delivered programs \citep{fukuoka2014b}.  The second type of incentive is that the clinician can provide a limited amount of counseling to individuals, but this is costly and the clinician must decide how to allocate a limited budget of counseling sessions to the entire pool of individuals.  For example, the \cite{dpp2002,dpp2003,dpp2009} has shown that such programs lead to a clinically significant loss of 5-7\% body weight on average, which can prevent or delay the onset of type 2 diabetes with few side effects.  However, these programs are difficult to design because variations in individual motivational states mean there is not just one set of optimal behavioral goals and assignment of counseling sessions, but rather that behavioral goals and the number/timing of counseling needs to be personalized to individuals' motivational states to maximize weight loss.

Personalizing the behavioral incentives for each individual can improve efficacy of weight loss programs and reduce the associated program costs through a reduction in the average amount of counseling for each individual.  Mobile phone technology is one promising avenue for implementing such personalization, due to its relatively low cost and pervasiveness among diverse communities \citep{lopez2013}. Mobile phones allow clinicians to collect real time health data through use of personal logs and devices such as  accelerometers, which provides noisy measurements of the health state and decisions of each individual.  Randomized controlled trials (RCT's) have found that the use of mobile phones can reduce the cost of implementing weight loss programs with maintaining efficacy \citep{fukuoka2014b}; however, little research to date has explored how to use the data generated by mobile phones and digital accelerometers in order to personalize behavioral incentives \citep{fukuoka2011,oreilly2013,azar2013,pagoto2013}. 
%In previous work, we proposed a predictive behavioral approach that utilized such data in order to make effective individual-specific predictions for weight loss using only a short time span of data \citep{aswani2016}. In this paper, we expand upon this methodology and show how it allows for the calculation of efficient and effective personalized interventions. Furthermore, we examine how these personalized interventions can be combined to form weight loss programs that maximize overall program effectiveness while maintaining cost efficiency.
 
Several adaptive methods have been proposed for designing personalized healthcare treatments, including: multi-armed bandits \citep{negoescu2014, deo2013,bastani2015online}, robust optimization \citep{ohair2013}, and dynamic programming \citep{engineer2009}. One common approach for optimal treatment design and clinical appointment scheduling has been Markov decision process (MDP) models \citep{ayer2015, mason2013, deo2013, kucukyazici2011,leff1986,liu2010,wang2011, gupta2008b,smilowitz2016}.  These methods are designed for situations with infrequent data collection (e.g., only during clinical visits), whereas in weight loss programs the data is collected daily (or more often) using mobile devices.  Our work develops an approach that can leverage this increased data availability to better design incentives.  Moreover, existing approaches focus either on motivational states characterizing adherence \citep{mason2013} or health states describing prognosis \citep{ayer2012,deo2013,helm2015,wu2013,negoescu2014,engineer2009}.  In contrast, we seek to combine motivational and health states into a single predictive model that is used for personalizing the weight loss program.  

\subsubsection{Demand Response (DR)}
DR programs are used by electric utilities to alter homeowners' electricity usage to better match electricity generation and reduce peak electricity demand.  Utilities incentivize homeowners to shift or reduce electricity consumption using price-based programs (e.g., time-of-day electricity rates).  Utilities also incentivize reduced electricity consumption through exchange programs in which a homeowner's inefficient appliances are replaced (for free by the electric utility) with efficient appliances \citep{palensky2011,deng2015}. Implementations of such DR programs have decreased peak electricity demand by almost 10\% and have improved the balance between electricity supply and demand \citep{lee2014}.  However, adverse selection is a major issue in these programs because incentives are often provided to homeowners who already have low electricity consumption or already had plans to replace inefficient appliances.

Better targeting in a DR program may lead to improved efficacy with lower associated costs. For example, electric utilities have the capability to send an auditor to homes to assess what appliance upgrades are needed \citep{website:PGE}. Consequently, an electric utility would be interested in finding the most effective way to schedule its auditors and set its rebates and tariffs.  Homeowner electricity usage data can be collected by the utility in real time using smart electricity meters, and the adoption rate of these smart electricity meters is increasing in the US \citep{lee2014}.  Moreover, the two way communication capabilities of smart electricity meters and mobile phones can be used to communicate billing and incentive information to homeowners \citep{Darby2010}, which opens the possibility for better targeting of price-based programs and appliance-replacement programs.

DR programs are often designed using game-theoretic approaches \citep{saghezchi2015,samadi2010,samadi2012}, multi-armed bandits \citep{wijaya2013}, convex optimization \citep{li2011,mohsenian2010,ratliff2014incentive}, dynamic programming \citep{jiang2011,costa2007,molderink2010}, and MDP's \citep{oneill2010,kim2011}.  These approaches commonly assume the electric utility has perfect information on the motivational state of each homeowner, and that the uncertainty is primarily in electricity generation and pricing.  In contrast, our proposed methodology has the ability to estimate the motivational state of each homeowner to better design DR programs through improved targeting of price-based and appliance-replacement incentives.  (Existing work also does not consider the option of the power company to provide rebates for upgrading inefficient equipment, while our framework can incorporate this scheduling problem.) 
  
\subsubsection{Franchise Logistics}
A franchise contract governs the relationship between a parent company (with a developed brand/products/services) and a separate franchisee company that retails this brand for a limited amount of time \citep{rubin1978}. Franchisees enter these contracts because they lack the brand recognition, economies of scale, or expertise needed.  Each franchisee desires to maximize their individual profit while the parent company works to ensure brand quality and profitability are maintained.  Towards this end, franchise contracts have several common clauses that help the parent company ensure that franchisees operate within quality standards: The parent company provides managerial assistance, training, and funding for advertising campaigns \citep{rubin1978}, and the parent company may also monitor a specific franchisee by assigning an employee to supervise the franchisee's daily operations \citep{gal1995}.

A common feature of such franchise contracts is the requirement for the franchisee to provide sales figures, customer satisfaction metrics, and expense reports to the parent company; and new online survey tools give the capability to monitor customer satisfaction and sales at a finer time-scale than previously possible.  Such data lends itself towards approaches for the improved allocation of parent company resources (e.g., managerial assistance, training, funding for advertising campaigns, and monitoring of specific franchisees), which can improve the performance of both the parent company and the franchisee.  Restated, these new data streams can be used by the parent company to develop a better understanding of the competencies and weaknesses of each franchisee, and to then use this to better allocate resources and set contract terms.

Characteristics of optimal franchise contracts have been studied extensively in the operations research literature \citep{gal1995,lal1990,li1997transaction,xie2016}. This stream of literature primarily extends earlier approaches used to solve contracting problems in which it is assumed that data acquisition and monitoring is expensive, making it difficult to estimate the utility function of the franchisee \citep{lal1990,gal1995}. Additionally, much of this work is focused on modeling the interaction and motivations of the franchisees and parent company and not on developing algorithms for optimal contract calculation.  In contrast, our behavioral analytics framework for adaptively designing incentives -- which consists of repeatedly estimating utility functions and then refining the incentives using optimization modeling -- applies to settings where data (e.g. sales, customer satisfaction, etc.) is readily available to the parent company.

\subsection{Literature Review}

The behavioral analytics framework we develop in this paper builds upon existing literature on data driven and adaptive methods for stochastic optimization. \cite{vahn2014} and \cite{vahn2015} consider how predictive and data driven models can be incorporated into inventory management problems, and both parametric and nonparametric predictive models are used by a decision maker to estimate demand and compute an optimal reorder policy.  These models are constructed to estimate demand through i.i.d observations; this differs from the setup in our paper where the observations are generated by temporal dynamics and are thus not i.i.d.  A more general set of approaches are reinforcement learning and Bayesian optimization \citep{aswani2013_automatica,frazier2016,osband2015,osband2016}, which leverage statistical estimation to compute asymptotically optimal control inputs for systems with appropriate model structures.  However, the relationship between the computed control inputs and the estimated model is often difficult to interpret because of the nonparametric nature of the estimation \citep{breiman2001}.  Our approach offers improved interpretability of the incentives computed by our framework because we simultaneously generate estimates of the parameters of the utility function (i.e., motivational states) for each agent.  These estimates provide insights into the resulting incentive allocations computed by our framework because these parameters usually have behavioral or financial interpretations (e.g., responsiveness to incentives, production efficiency, level of risk aversion).

Our behavioral analytics framework is also related to research that explores stochastic control of multi-agent systems. Related methods include decentralized control \citep{li2012}, approximate dynamic programming \citep{boukhtouta2011,george2007}, game-theoretic approaches \citep{adlakha2013,iyer2011,iyer2014,Zhou2016}, and robust optimization \citep{blanchet2013,bertsimas2012,lorca2015}. In general, these models consider very different settings from the ones we consider in this paper.  This body of work studies settings where the agents can strategically interact with other agents (without the presence of a coordinator) and where the agents are able to consider long time horizons when making decisions.  Our setting differs in that we have a single coordinator that provides incentives to a group of agents that do not interact strategically with other agents, and where the agents are myopic (meaning they make decisions based on short time horizons).  The three examples of weight loss programs, demand response programs, and franchise logistics more closely match the setting we consider in this paper.

\subsection{Contributions}

Our overall goal in this paper is to provide tools and approaches that form a specific implementation of the three steps of a behavioral analytics framework, and our secondary goal is to give an example that demonstrates how our implementation of behavioral analytics can be applied to a real-world engineering problem.  Recall that these three steps involve designing a behavioral model, and then repeatedly estimating the parameters of this behavioral model, and using the estimated parameters to optimize the incentives provided to each agent.  To do this, we first need to identify a general (and practically useful) class of models that describe agent behavior and can be incorporated into optimization models for incentive design.  This is non-obvious because incentive design in principle requires solving bilevel programs, precluding the straightforward use of commercial optimization software packages.  We address this by abstracting and generalizing our earlier work on the development of predictive models of the behavior of individuals participating in a weight loss program \citep{aswani2016}. Given these behavioral models, we design an optimization approach that, rather than directly solving the relevant bilevel program,  is built around formulations that incorporate the individual behavior model into mathematical programs that can be solved in a straightforward way with commercial solvers, and that lead to incentives that are asymptotically optimal as more data is collected.  Below, we describe these contributions in further detail:

First, we develop and analyze an abstract model of agent behavior.  This model consists of a myopic utility function (meaning the agent makes decisions based on a utility function that depends on states only one time period into the future) and temporal dynamics on the system states and on the parameters of the utility function.  It abstracts and generalizes a predictive model we created in our prior work on behavioral modeling for weight loss \citep{aswani2016}.  In addition, we explore (for the first time) theoretical questions related to statistical consistency of utility function parameter estimates.  Such consistency is important because in order to design optimal incentives we need to be able to correctly estimate the parameters of the utility functions of each agent, and it was recently shown that not all approaches that have been proposed for estimating parameters of utility functions are statistically consistent \citep{siddiq2015}.  Here, we provide mixed integer linear programming (MILP) formulations for estimating the parameters of the utility functions, and we prove these formulations generate estimates that are statistically consistent.

% This paper also makes contributions towards solving bilevel mixed-integer programs (BMIP's) through development of a two-stage decomposition algorithm for BMIP's with a specific mathematical structure that arises from the coordinator's problem because the distributions of future agent states are represented by nonparametric Bayesian distributions computed by MILP's.  BMIP's are computationally difficult to solve \citep{ralphs2014,denegre2009,moore1990,moore1992} since solution techniques for continuous bilevel programming \citep{ahuja2001,siddiq2015,aswani2016duality,dempe2002,heuberger2004} cannot be used.  Our two-stage decomposition algorithm first solves as sub-problems the coordinator's problem for each agent considered individually, and we show these sub-problems can be posed as MILP's.  Next, our two-stage decomposition algorithm solves a single master problem that can be cast as a pure integer linear program (ILP).  We prove that if the Bayesian distributions can be described by a sufficient and consistent statistic, then the overall algorithm generates solutions that are asymptotically optimal as time $t$ goes to infinity.

We also develop novel mathematical programs for incentive design that incorporate our model for agent behavior, and we prove that the incentives are asymptotically optimal (in time).  Incentive design in principle requires solving a bilevel program, and the situation is complicated in our setting because the mathematical structure of our abstract model for agent behavior leads to incentive design problems that consist of bilevel mixed integer programs (BMIP's). BMIP's are computationally difficult to solve \citep{ralphs2014,denegre2009,moore1990,moore1992} since solution techniques for continuous bilevel programming \citep{ahuja2001,siddiq2015,aswani2016duality,dempe2002,heuberger2004} cannot be used. Consequently, we develop an adaptive two-stage decomposition algorithm.  In the first stage, we solve the coordinator's problem for each agent considered individually by estimating the utility function parameters of an agent by solving a single MILP and then solving a series of MILP sub-problems.  The second stage consists of an integer linear program (ILP) master problem that aggregates the sub-problem solutions and solves the coordinator's problem for all agents considered jointly.  We prove this asymptotically designs the optimal incentives.

To evaluate the efficacy of the specific behavioral models, parameter estimation techniques, and optimization models in our instantiation of a behavioral analytics framework, we perform computational experiments in the context of goal-setting and clinical appointment scheduling for individuals participating in a clinically-supervised weight loss program.  The first step of our behavioral analytics approach involves constructing a model that describes individuals' decisions on how much to eat and how much physical activity (in terms of daily steps) to do -- subject to a utility function that captures the tradeoffs inherent in achieving one-day-ahead weight loss with reducing dietary consumption and increasing physical activity.  The second step of our behavioral analytics approach uses past data for each individual in order to quantify (for each individual) the tradeoffs captured by the utility function, as well as estimate the responsiveness of each individual to the incentives of providing physical activity goals and providing counseling sessions, and the third step of our behavioral analytics approach uses the behavioral model and estimated parameters to determine what physical activity goals to provide to each individual and to determine how to allocate a fixed number of counseling sessions to a pool of individuals participating in the program.  These second and third steps are repeated as more data is collected from each individual. Through a simulation study, we compare personalized treatment plans computed by our approach with treatment plans computed by an adaptive heuristic, and we find that our approach performs substantially better than the heuristic.  Common heuristics implicitly assume monotonicity in individuals' behaviors with respect to the treatment plan values, while actual behavior (captured by our predictive models) displays substantial non-monotonicity: For example, losing weight causes individuals to eat more and exercise less, so the speed of weight loss can impact the final weight loss outcomes.

 \subsection{Outline}

%Section \ref{sec:mult_aget} presents our two-stage algorithm for solving the coordinator's problem, and Sections \ref{sec:pred_mod} and \ref{sec:opt_sing_age} provide the elements required to describe the two-stage algorithm.  

Section \ref{sec:pred_mod} describes the first step of our behavioral analytics framework -- the development of the behavioral model.  The model consists of a utility function -- describing how an agent makes decisions -- and temporal dynamics on the system states and parameters of the utility function.  We refer jointly to both components of this abstract model as the behavioral model.  Section \ref{sec:pred_mod_est} presents approaches for estimating parameters of this behavioral model using MILP formulations to solve the problems of maximum likelihood estimation (MLE) and Bayesian inference.  We prove that solutions of our MILP formulations provide consistent estimates of the agent's parameters.  In Section \ref{sec:opt_sing_age}, we present algorithms for optimizing the incentives provided to agents by the coordinator.  We first present an algorithm based on solving two MILP's that allows the coordinator to allocate incentives in the situation where there is only a single agent with unknown-to-the-coordinator parameters, and we prove that this algorithm computes incentives that are asymptotically optimal (in the sense of minimizing the coordinator's loss function) as time $t$ goes to infinity.  Next, we develop a two-stage decomposition algorithm (building on the single-agent formulation) to solve the coordinator's problem in a multi-agent setting, and we generalize our proof of asymptotic optimality to this setting  Finally in Section \ref{sec:experiment}, we study (via simulation) the effectiveness of our algorithms for designing personalized weight loss treatment plans. Our results show that treatment plans computed by our behavioral analytics approach could potentially reduce the cost of running such weight loss programs by as much as 60\% without affecting the efficacy of these programs.
 
\section{Predictive Modeling of a Single Myopic Agent}
\label{sec:pred_mod}
In this section, we present our behavioral model for a single myopic agent.  This forms the first step of our specific implementation of a behavioral analytics framework, and the key design problem is formulating a predictive model that is amenable to performing the second and third steps of our behavioral analytics framework of parameter estimation and incentive optimization.   This model is an extension and abstraction of a behavioral model that was validated in our past work on behavioral modeling for weight loss \citep{aswani2016}, in which we used cross-validation (i.e., out-of-sample comparisons) to perform a data-based validation of the predictive accuracy of our behavioral model by comparison to a standard machine learning algorithm for prediction.  %The behavioral model presented here abstracts the model by \cite{aswani2016}; we also prove theoretical results on statistical consistency in this paper, which was not considered in our previous work \citep{aswani2016}.

%\subsection{Behavioral Model}
Let $\mathcal{X},\mathcal{U},\Pi,\Theta$ be compact finite-dimensional sets with $\mathcal{X},\mathcal{U},\Theta$ convex. We will refer to the agent's system states $x_t\in\mathcal{X}$, motivational states (or type) $\theta_t\in\Theta$, and decisions $u_t \in \mathcal{U}$ at time $t$. The coordinator provides an incentive (or input) $\pi_t \in \Pi$ to the agent at time $t$, and we assume that the motivational states are unknown to the coordinator but known to the agent.  In our behavioral model, the system and motivational states are subject to temporal dynamics:
\begin{equation} 
\label{eq:dynam_eq}
\begin{aligned}
& x_{t+1} = h(x_t,u_t), \\
&\theta_{t+1} = g(x_t,u_t,\theta_t,\pi_t).
\end{aligned}
\end{equation}
The intuition of the above dynamics is that future system states $x_{t+1}$ depend on the current system states $x_t$ and decision $u_t$, while future motivational states $\theta_{t+1}$ depend on the current system states $x_t$, decision $u_t$, motivational states $\theta_t$, and incentives $\pi_t$.  

The agents are modeled to be myopic in the sense that agents make decisions at time $t$ by considering only their present utility function.  We assume the agent's utility function belongs to a parametrized class of functions $\mathcal{F}:=\{(x,u) \mapsto f(x,u,\theta,\pi):\theta\in\Theta, \pi \in \Pi\}$; and the agent's utility function at time $t$ is $f(\cdot,\cdot,\theta_t,\pi_t)$. Thus at time $t$ the agent's decisions are
\begin{equation}
\label{eq:dec_eq}
u_t \in \argmax\big\{f(x_{t+1},u,\theta_t,\pi_t)\ \big|\ x_{t+1} = h(x_t,u),\ u\in\mathcal{U}\big\},
\end{equation}
which means we are assuming the agent has perfect knowledge of $x_t,\theta_t,\pi_t$.  This model says that the agent's decisions depend on the current system states, motivational states, and incentives.  

Though the coordinator also has perfect knowledge of the incentives $\pi_t$, the coordinator can only make noisy observations of past system states and agent decisions:
\begin{equation}
\begin{aligned}
\tilde{x}_{t_i} &= Dx_{t_i} + \nu_{t_i} &\forall i = 0,\ldots,n_x,\\
\tilde{y}_{\tau_i} &= Cu_{\tau_i} + \omega_{\tau_i} &\forall i = 0,\ldots,n_u,
\end{aligned}
\label{eq:mle_noise}
\end{equation}
where $C,D$ are known output matrices, and $x_{t_i},u_{\tau_i}$ are the systems states and agent decisions generated by (\ref{eq:dynam_eq}) and (\ref{eq:dec_eq}) with initial conditions $(x_0,\theta_0)$ and incentives $(\pi_1,\pi_2,\ldots)$. 

For our subsequent optimization modeling and theoretical analysis, we make the following assumptions about this behavioral model:
%\begin{assumption}
%$\forall (x,\theta,\pi) \in \mathcal{X} \times \Theta \times \Pi$, the set $\argmax_{u\in\mathcal{U}}f(x,u;\theta,\pi)$ is singleton, \label{ass:ass_1}
%\end{assumption}
\begin{assumption}
\label{ass:ass_1}
The sets $\mathcal{X},\mathcal{U},\Pi,\Theta$ are bounded and finite-dimensional. Moreover, the sets $\mathcal{X},\mathcal{U},\Theta$ are convex polyhedra described by a finite number of linear inequalities, and $\Pi$ can be described by a finite number of mixed integer linear constraints.
\end{assumption}
This mild assumption ensures that states, decisions, and inputs are bounded; that the range of possible values for states and inputs are polytopes; and that the set of possible incentives is representable by mixed integer linear constraints.
\begin{assumption}
\label{ass:ass_2}
%The function $f$ is deterministic, concave in $x$, strictly concave in $u$, concave in $\theta$, and dependent on $\pi,\theta$ through linear and bilinear terms; moreover, $f$ can be expressed as a sum of quadratic terms and continuous piecewise linear terms;
The function $f : \mathcal{X}\times\mathcal{U}\times\Theta\times\Pi\rightarrow\mathbb{R}$ is deterministic, concave in $x$, strictly concave in $u$, and concave in $\theta$; moreover, $f$ can be expressed as 
\begin{equation}
f(x,u, \theta, \pi) = -(x; u)^T\cdot Q\cdot(x;u) + (\theta;\pi)^T\cdot H\cdot(x;u) + \sum_{i=1}^K \min_{j\in J_i}\big\{F_{i,j}\cdot(x;u;\theta;\pi) + \zeta_{i,j}\big\},
\end{equation}
where $Q$ is a positive semidefinite matrix, the $F_{i,j},H$ are matrices of appropriate dimension, the $\zeta_{i,j}$ are scalars, and the $J_i$ are sets of indices.
\end{assumption}
Strict concavity in $u$ ensures $\argmax_{u\in\mathcal{U}}f(x,u,\theta,\pi)$ is singleton for all $(x,\theta,\pi) \in \mathcal{X} \times \Theta \times \Pi$, and the concavity assumptions also model diminishing returns and ensure $u_t$ is polynomial-time computable by the agent \citep{brock1990,gafni1990,cawley2004}.%  The remaining characteristics of $f$ allow us to formulate problems of statistical estimation as a MILP.
\begin{assumption}
\label{ass:ass_3}
The functions $h : \mathcal{X}\times\mathcal{U}\times\Theta\times\Pi\rightarrow\mathcal{X}$ and $g : \mathcal{X}\times\mathcal{U}\times\Theta\times\Pi\rightarrow\Theta$ are deterministic surjective functions of the form
\begin{equation}
\begin{aligned}
&h(x,u) = Ax + Bu + k\\
&g(x,u,\theta,\pi) = G_i\cdot(x;u;\theta;\pi) + \chi_i \mathrm{\ when\ } B_i\cdot(x;u;\theta;\pi) \leq \psi_i
\end{aligned}
\end{equation}
where $A,B,G_i,B_i$ are matrices; $\gamma_i, \psi_i,k$ are vectors; $\chi_i$ are scalars; and the interiors of the polytopes $B_i\cdot(x;u;\theta;\pi) \leq \psi_i$ are disjoint.
\label{ass:ass_4}
\end{assumption}
This condition on $h,g$ allows us to formulate problems of statistical estimation as a MILP.
\begin{assumption}
The $\{\nu_{t_i}\}_{i=0}^{n_x}$ and $\{\omega_{\tau_i}\}_{i=0}^{n_u}$ from the measurement noise model \eqref{eq:mle_noise} are sequences of i.i.d random vectors with i.i.d components with zero mean and (known) finite variance.  Moreover, the logarithm of their probability density functions can be expressed using integer linear constraints.
\end{assumption}
This means $\mathbb{E}\omega_{\tau_i} = \mathbb{E} \nu_{t_i} =0$ and $\mathbb{E}(\nu_{t_i})_j^2 = \sigma_\nu^2 < \infty$ and $\mathbb{E}(\omega_{\tau_i})_j^2 = \sigma_\omega^2 < \infty$ with known $\sigma_\nu^2,\sigma_\omega^2$.  Examples of noise distributions satisfying the integer linear representability assumption include the Laplace distribution, the shifted exponential distribution, and piecewise linear distributions.  This assumption can be relaxed to requiring integer quadratic representability (such as is the case for Gaussian distributed noise), and the subsequent results change in that the optimization formulations become MIQP's, rather than the MILP's that occur with the above assumption.
%\begin{assumption}
%For each fixed $\pi_t\in\Pi$ the mapping from $(x_t,\theta_t)$ to $(x_{t+1},\theta_{t+1})$ defined by \eqref{eq:dynam_eq} and \eqref{eq:dec_eq} is bijective.
%\end{assumption}
%This assumption implies that each distinct initial condition $(x_0,\theta_0)$ uniquely defines a future trajectory of system states $x_t$ and motivational states $\theta_t$ when the incentives $(\pi_1,\pi_2,\ldots)$ are fixed.
\begin{assumption}
The discrete-time system with temporal dynamics \eqref{eq:dynam_eq} and \eqref{eq:dec_eq} and measurement model \eqref{eq:mle_noise} is observable (i.e., there exists a $T$ and sequence $\pi_t$ such that $(x_0,\theta_0)$ can be exactly computed if the measurements from $0 \leq t \leq T$ are noiseless).
\label{ass:ass_5} 
\end{assumption}
This last assumption is an identifiability condition \citep{bickel2006}, meaning that different initial conditions $(x_0,\theta_0)$ on the agent's system and motivational states produce different sequences of measurements and states, and this is a common assumption for control systems \citep{callier1994}.  The second assumption is common for utility functions \citep{brock1990,gafni1990,cawley2004}.  The third assumption says the system state dynamics are linear, and that the motivational state dynamics are piecewise affine, which are common models for control systems \citep{callier1994,mignone2000,aswani2009monotone}.  We believe all five assumptions are satisfied by agents in the three examples of weight loss programs, demand response programs, and franchise logistics.  Section \ref{sec:experiment} provides a behavioral model for agents in a weight loss program that satisfies our above assumptions, and we conclude that section with a computational study where we solve the coordinator's problem for a weight loss program.

% More broadly, our assumptions imply that the initial conditions are sufficient to describe the entire future progression of the agent's behavior for a given sequence of incentives $(\pi_1,\pi_2,\ldots)$.  

\section{Estimating Model Parameters}
\label{sec:pred_mod_est}

In this section, we explore how the coordinator can estimate the agent's initial parameters $(x_0,\theta_0)$, and predict the agent's future behavior for a fixed policy $\pi$.  This forms the second step of our implementation of a behavioral analytics framework, and we leverage the mathematical structure of the behavioral model described in Section \ref{sec:pred_mod} to construct techniques and methods for estimation and prediction.  This second step is important because the estimated parameters of the behavioral model and subsequent predictions of future agent behavior are used to optimize incentives in the third step of our behavioral analytics approach. We will assume the coordinator makes noisy and partial observations -- according to the measurement model (\ref{eq:mle_noise}) -- of the agent's state and decisions for $n$ time periods (with some missing observations).   In Section \ref{sec:mle}, we present a Maximum Likelihood Estimation approach to  estimate the agent's initial system states and motivational states.  In Section \ref{sec:bayes_pred}, we consider a setting in which the coordinator has some prior knowledge about the possible values of the motivational states, and consider a Bayesian setting.  

\subsection{Maximum Likelihood Estimation}
\label{sec:mle}

Let $\{\tilde{x}_{t_i}\}_{i=0}^{n_x}$ denote the process of the state observations, and let $\{\tilde{y}_{\tau_i}\}_{i=0}^{n_u}$ denote the process of the behavior observations. 

Our approach to estimating the agent's initial parameters will be to compute estimates $(\hat{x}_0,\hat{\theta}_0) \in \argmin_{(x_0,\theta_0) \in \mathcal{X}\times\Theta} \mathcal{L}(x_0,\theta_0,\{\tilde{x}_{t_i}\}_{i=0}^{n_x},\{\tilde{y}_{\tau_i}\}_{i=0}^{n_u};\pi)$ by minimizing an appropriately chosen loss function $\mathcal{L}$.  More specifically, we use the approach of Maximum Likelihood Estimation (MLE).  Let $p_\nu, p_\omega$ be the density functions of $\nu_{t_i},\omega_{\tau_i}$; then the joint likelihood function of $(\theta_0, x_0)$ for a fixed $\pi$ is
\begin{equation}
\begin{aligned}
\mathcal{L}(x_0,\theta_0,\{\tilde{x}_{t_i}\}_{i=0}^{n_x},\{\tilde{y}_{\tau_i}\}_{i=0}^{n_u},\pi) &= p(\{\tilde{x}_{t_i}\}_{i=0}^{n_x},\{\tilde{y}_{\tau_i}\}_{i=0}^{n_u} | \theta_0,x_0,\pi)\\
 &= \prod_{i=0}^{n_x} p_\nu(\tilde{x}_{t_i}-Dx_{t_i}) \prod_{j=0}^{n_u} p_\omega(\tilde{y}_{\tau_j}-Cu_{\tau_j})
\end{aligned}
\end{equation}
Thus the coordinator's estimation problem is given by the following:
\begin{equation}
\begin{aligned}
(\hat{x}_0,\hat{\theta}_0) \in \argmax_{(x_t,\theta_t,u_t)}& \sum_{i=0}^{n_x} \log p_\nu(\tilde{x}_{t_i}-Dx_{t_i}) + \sum_{j=0}^{n_u} \log p_\omega(\tilde{y}_{\tau_j}-Cu_{\tau_j}) \\
 \text{s.t. } & u_t \in \argmax\big\{f(x_{t+1},u,\theta_t,\pi_t)\ \big|\ x_{t+1} = h(x_t,u),\ u\in\mathcal{U}\big\} & \forall t, \\
& x_{t+1} = h(x_t,u_t) & \forall t, \\
&\theta_{t+1} = g(x_t,u_t,\theta_t,\pi_t) &\forall t, \\
& x_t \in \mathcal{X}, \theta_t \in \Theta &\forall t.
\end{aligned}
\label{eq:mle_prob}
\end{equation}
Problem (\ref{eq:mle_prob}) is a bilevel optimization problem because the $u_t$ are minimizers of $f(x_{t+1},\cdot,\theta_t,\pi_t)$, and such bilevel problems frequently arise in the context of estimating utility functions \citep{keshavarz2011,bertsimas2014,siddiq2015}.  For the setting we consider in this paper, we show that the bilevel program for MLE \eqref{eq:mle_prob} can be exactly reformulated as a MILP.
\begin{proposition} \label{prop:form_const}
If Assumptions \ref{ass:ass_1}--\ref{ass:ass_5} hold; then the feasible region of \eqref{eq:mle_prob} can be formulated as a set of mixed integer linear constraints with respect to $(x_t,u_t,\theta_t,\pi_t)$. 
\end{proposition}

\proof{Proof: }
The constraints $\theta_{t+1} = g(x_t,u_t,\theta_t,\pi_t)$ can be reformulated using Assumption \ref{ass:ass_3} as
\begin{equation}
\begin{aligned}
\theta_{t+1} \leq G_i\cdot(x_t;u_t;\theta_t;\pi_t)+\xi_i + (1-\iota_i)\cdot M\\
\theta_{t+1} \geq G_i\cdot(x_t;u_t;\theta_t;\pi_t)+\xi_i - (1-\iota_i)\cdot M\\
B_i\cdot(x_t;u_t;\theta_t;\pi_t)\leq\psi_i + (1-\iota_i)\cdot M\\
\iota_i\in\{0,1\}
\end{aligned}
\end{equation}
where $M > 0$ is a large-enough constant.  Such a finite $M$ exists because $\mathcal{X},\mathcal{U},\Pi,\Theta$ are compact.  Hence it suffices to show $u_t \in \argmax\big\{f(x_{t+1},u,\theta_t,\pi_t)\ \big|\ x_{t+1} = h(x_t,u),\ u\in\mathcal{U}\big\}$ can be represented (by its optimality condition) using a finite number of mixed integer linear constraints. Suppose $\mathcal{U} = \{u : \Xi u \leq \kappa\}$, where $\Xi$ is a matrix and $\kappa$ is a vector.  Recall Assumption \ref{ass:ass_2}
\begin{equation}
\label{eqn:fst}
f(x,u, \theta, \pi) = -(x; u)^T\cdot Q\cdot(x;u) + (\theta;\pi)^T\cdot H\cdot(x;u) + \sum_{i=1}^K \min_{j\in J_i}\big\{F_{i,j}\cdot(x;u;\theta;\pi) + \zeta_{i,j}\big\}.
\end{equation}
We cannot characterize optimality by differentiating $f$ because it is generally not differentiable, but we can reformulate the maximization of \eqref{eqn:fst} as the following convex quadratic program:
\begin{equation}
\begin{aligned}
u_t \in \arg\max\ & -(x_{t+1}; u)^T\cdot Q\cdot(x_{t+1};u) + (\theta_t;\pi_t)^T\cdot H\cdot(x_{t+1};u) + \sum_{i=1}^K w_i \\
\text{s.t. }& w_i \leq F_{i,j}\cdot(x_{t+1};u;\theta_t;\pi_t) + \zeta_{i,j} \text{ for all } i,j\\
&\Xi u\leq\kappa 
\end{aligned}
\end{equation}
Using Assumption \ref{ass:ass_4}, we can rewrite the above as
\begin{equation}
\label{eqn:reform2}
\begin{aligned}
u_t \in \arg\max\ & -u^T\cdot(B; \mathbb{I})^T\cdot Q\cdot(B; \mathbb{I})\cdot u + ((\theta_t;\pi_t)^T H-2(Ax+k; 0)^T Q)\cdot(B; \mathbb{I})\cdot u+ \sum_{i=1}^K w_i \\
%\text{s.t. }& w_i \leq (\theta_t;\pi_t)^T\cdot E_{i,j}\cdot(Ax_t + Bu +k;u) + F_{i,j}^T\cdot(Ax_t+Bu+k;u) + \zeta_{i,j} \text{ for all } i,j 
\text{s.t. }& w_i \leq F_{i,j}\cdot(B; \mathbb{I}; 0; 0)\cdot u +  F_{i,j}\cdot(Ax_t+k; 0;\theta_t;\pi_t)+\zeta_{i,j} \text{ for all } i,j\\
 &\Xi u\leq\kappa 
\end{aligned}
\end{equation}
where we have eliminated the constant $(\theta_t;\pi_t)^T\cdot H\cdot(Ax_{t}+k;0)$ since $x_t,\theta_t,\pi_t$ are known to the agent. The above optimization problem is strictly convex by assumption, and all constraints are linear for fixed $\pi_t$.  Hence the optimality conditions for (\ref{eqn:reform2}) can be characterized using the KKT conditions \citep{dempe2002,boyd2004}. Let $\lambda_{i,j}$ and $\mu$ be the Lagrange Multipliers for the first and second set of constraints given in \eqref{eqn:reform2}, and note the KKT conditions are
\begin{equation}
\begin{aligned}
& 2(B; \mathbb{I})^T\cdot Q\cdot(B; \mathbb{I})\cdot u + (B; \mathbb{I})^T\cdot (2Q\cdot(Ax+k; 0)-H^T\cdot(\theta_t;\pi_t)) + \Xi^T\cdot \mu = \\
&\hspace{8.5cm}\sum_{i=1}^K\sum_{j\in J_i}\lambda_{i,j}\cdot(B; \mathbb{I};0;0)^T\cdot F_{i,j}^T\\
&\sum\limits_{j\in J_i} \lambda_{i,j} = 1 \text{ for } i = 1, \ldots,K \\
&\lambda_{i,j}\cdot(w_i - F_{i,j}\cdot(B; \mathbb{I}; 0; 0)\cdot u +  F_{i,j}\cdot(Ax_t+k; 0;\theta_t;\pi_t)+\zeta_{i,j}) = 0 \text{ for all } i,j  \\
 &w_i \leq F_{i,j}\cdot(B; \mathbb{I}; 0; 0)\cdot u +  F_{i,j}\cdot(Ax_t+k; 0;\theta_t;\pi_t)+\zeta_{i,j} \text{ for all } i,j\\
 &\Xi u\leq\kappa\text{ and }\lambda_{i,j} \geq 0 \text{ for all } i,j
\end{aligned}
\end{equation}
Note that the only nonlinear conditions are those which represent complimentary slackness. However, these conditions can be reformulated as integer linear conditions by posing them as disjunctive constraints \citep{wolsey1999}: For sufficiently large $M$ -- which exists because of the compactness of $\mathcal{X},\mathcal{U},\Pi,\Theta$  -- the complimentary slackness conditions are
\begin{equation}
\begin{aligned}
&\lambda_{i,j} \leq M \iota_{i,j} \text{ for all } i,j \\
&F_{i,j}\cdot(B; \mathbb{I}; 0; 0)\cdot u +  F_{i,j}\cdot(Ax_t+k; 0;\theta_t;\pi_t)+\zeta_{i,j} \leq w_i + M\cdot(1-\iota_{i,j}) \text{ for all } i,j\\
&\iota_{i,j} \in \{0,1\} \text{ for all } i,j
\end{aligned}
\end{equation}
This shows that feasible region of \eqref{eq:mle_prob} of can be represented using a finite number of mixed integer linear constraints.\Halmos
\endproof
An important consequence of this proposition is that it is possible to compute the global solution of the MLE problem (\ref{eq:mle_prob}) using standard optimization software.

\begin{corollary}\label{cor:mle_cor1}
If Assumptions \ref{ass:ass_1}--\ref{ass:ass_5} hold, then the MLE problem \eqref{eq:mle_prob} can be expressed as a MILP.
\end{corollary}

\begin{remark}
If the logarithm of the noise densities can be expressed using integer quadratic constraints (e.g., Gaussian distributions), then the MLE problem \eqref{eq:mle_prob} can be expressed as a MIQP.
\end{remark}

\subsection{Bayesian Estimation}
\label{sec:bayes_pred}
Solving the MLE problem \eqref{eq:mle_prob} gives an estimate of the agent's initial system states and motivational states, which completely characterize the agent.  However, the coordinator often has some prior knowledge about the possible values of the motivational states.  In such a case, a Bayesian framework is a natural setting for making predictions of the agent's future system states.  

Suppose the coordinator has interacted with the agent over $T$ time periods, has measured $\{\tilde{x}_{t_i}\}_{i=0}^{n_x}, \{\tilde{y}_{\tau_i}\}_{i=0}^{n_u}$ with $T_{n_x},T_{n_u} \leq T$, and wants to predict the agent's future states and decisions $\{x_i,\theta_i,u_i\}_{i=T}^{T+n}$ for some $n > 0$ time steps into the future.  In principle, this means the coordinator wants to calculate the posterior distribution of $\{x_i,\theta_i,u_i\}_{i=T}^{T+n}$. But $(x_0,\theta_0)$ completely characterize the agent in our model (recall Assumption \ref{ass:ass_5} states that distinct initial conditions produce different state and decision trajectories), and so we can predict the agent's future states and decisions using the posterior distribution of $(x_0,\theta_0)$.  Hence we focus on computing the posterior of $(x_0,\theta_0)$.  A direct application of Bayes's Theorem \citep{bickel2006} gives
\begin{equation}
\label{eqn:bigpost}
p\big(x_0,\theta_0\big|\{\tilde{x}_{t_i}\}_{i=0}^{n_x}, \{\tilde{y}_{\tau_i}\}_{i=0}^{n_u},\{\pi_i\}_{i=0}^{T+n}\big) = Z^{-1}\times p\big(\{\tilde{x}_{t_i}\}_{i=0}^{n_x}, \{\tilde{y}_{\tau_i}\}_{i=0}^{n_u} \big| x_0,\theta_0,\{\pi_i\}_{i=0}^{T+n}\big)\times p\big(x_0,\theta_0\big).
\end{equation}
Here $Z$ is a normalization constant that ensures the right hand side is a probability distribution, and $p(x_0,\theta_0)$ reflects the coordinator's prior beliefs.  We begin with an assumption on $p(x_0,\theta_0)$.
\begin{assumption}
\label{ass:ass_6}
The function $\log p(x_0,\theta_0)$ can be expressed using a finite number of mixed integer linear constraints, and $p(x_0,\theta_0)>0$ for all $(x_0,\theta_0)\in \mathcal{X}\times\Theta$.
\end{assumption}
This is a mild assumption because it holds for the Laplace distribution, the shifted exponential distribution, and piecewise linear distributions.  Significantly, it is true when the prior distribution $p(x_0,\theta_0)$ is an empirical histogram with data in each histogram bin \citep{aswani2016}.

Next, we describe an optimization approach to computing the posterior distribution of $(x_0,\theta_0)$.  Consider the following feasibility problem:
\begin{equation}
\begin{aligned}
\psi_T(\bar{x}_0,\bar{\theta}_0) &= \log p\big(x_0=\bar{x}_0,\theta_0=\bar{\theta}_0\big|\{\tilde{x}_{t_i}\}_{i=0}^{n_x}, \{\tilde{y}_{\tau_i}\}_{i=0}^{n_u},\{\pi_i\}_{i=0}^T\big) + \log Z\\
&= \max_{(x_t,\theta_t,u_t)}\ \sum_{i=0}^{n_x} \log p_\nu(\tilde{x}_{t_i}-Dx_{t_i}) + \sum_{j=0}^{n_u} \log p_\omega(\tilde{y}_{\tau_j}-Cu_{\tau_j}) + \log p(x_0,\theta_0) \\
 &\hspace{1cm}\text{s.t. }\  \ u_t \in \argmax\big\{f(x_{t+1},u,\theta_t,\pi_t)\ \big|\ x_{t+1} = h(x_t,u),\ u\in\mathcal{U}\big\} & \forall t, \\
&\hspace{1.9cm}x_{t+1} = h(x_t,u_t) & \forall t, \\
&\hspace{1.9cm}\theta_{t+1} = g(x_t,u_t,\theta_t,\pi_t) & \forall t, \\
&\hspace{1.9cm}x_0=\bar{x}_0,\theta_0=\bar{\theta}_0, \\
&\hspace{1.9cm}x_t \in \mathcal{X}, \theta_t \in \Theta &\forall t.
\end{aligned}
\label{eq:prof_lik_prob}
\end{equation}   
The above problem is almost the same as the MLE problem (\ref{eq:mle_prob}), with the only differences that the above has additional constraints $x_0=\bar{x}_0,\theta_0=\bar{\theta}_0$ and an additional term in the objective $\log p(x_0,\theta_0)$.  Thus we have that the above problem (\ref{eq:prof_lik_prob}) can be expressed as a MILP or MIQP.
\begin{corollary}
If Assumptions \ref{ass:ass_1}--\ref{ass:ass_6} hold, then \eqref{eq:prof_lik_prob} can be formulated as a MILP.
\label{cor:mapmilp}
\end{corollary}

\begin{remark}
Under appropriate relaxed representability conditions on the noise distributions and the prior distribution, the problem \eqref{eq:prof_lik_prob} can be formulated as a MIQP.
\end{remark}

Solving \eqref{eq:prof_lik_prob} does not directly provide the posterior distribution of $(x_0,\theta_0)$ because $Z$ is not known \emph{a priori}, though it can be computed using numerical integration.  (See for instance the approach by \cite{aswani2016}.)  But since $Z$ only scales the posterior estimate, we instead propose a simpler scaling.  Let $(\hat{x}_{0,T},\hat{\theta}_{0,T}) \in \argmax_{(x_0,\theta_0)} \psi_T(x_0,\theta_0)$ be the maximum \emph{a posteriori} (MAP) estimates of the initial conditions, and note that the above corollaries apply to the computation of the MAP because the corresponding optimization problem for computing the MAP is simply (\ref{eq:prof_lik_prob}) but with the constraints $x_0=\bar{x}_0,\theta_0=\bar{\theta}_0$ removed.  We propose using
\begin{equation}
\hat{p}\big(x_0,\theta_0 \big|\{\tilde{x}_{t_i}\}_{i=0}^{n_x}, \{\tilde{y}_{\tau_i}\}_{i=0}^{n_u},\{\pi_i\}_{i=0}^T \big)= \frac{\exp(\psi_T(x_0,\theta_0))}{\exp(\psi_T(\hat{x}_{0,T},\hat{\theta}_{0,T}))} \label{eq:post_dens}
\end{equation}
as an estimate of the posterior distribution of $(x_0,\theta_0)$.  Two useful properties of our estimate are that $\hat{p}(x_0,\theta_0 |\{\tilde{x}_{t_i}\}_{i=0}^{n_x}, \{\tilde{y}_{\tau_i}\}_{i=0}^{n_u},\{\pi_i\}_{i=0}^T) \in[0,1]$ by construction, and that $\hat{p}(\hat{x}_{0,T},\hat{\theta}_{0,T}|\{\tilde{x}_{t_i}\}_{i=0}^{n_x}, \{\tilde{y}_{\tau_i}\}_{i=0}^{n_u},\{\pi_i\}_{i=0}^T)=1$ by construction.  We will show this estimate is statistically consistent in a Bayesian sense \citep{bickel2006}:
\begin{definition}
The posterior estimate (\ref{eq:post_dens}) is consistent if for all $(x_0^*,\theta_0^*) \in \mathcal{X}\times\Theta$ and $\epsilon,\delta>0$ we have $p_{(x_0^*,\theta_0^*)}(\hat{p}(\mathcal{E}(\delta)|\{\tilde{x}_{t_i}\}_{i=0}^{n_x}, \{\tilde{y}_{\tau_i}\}_{i=0}^{n_u},\{\pi_i\}_{i=0}^T) \geq \epsilon) \rightarrow 0$ as $T\rightarrow\infty$, where $p_{(x_0^*,\theta_0^*)}$ is the probability law under $(x_0^*,\theta_0^*)$, $\mathcal{E}(\delta) = \{(x_0,\theta_0) \notin \mathcal{B}(x_0^*,\theta_0^*,\delta)\}$, and $\mathcal{B}(x_0^*,\theta_0^*,\delta)$ is an open $\delta$ ball around $(x_0^*,\theta_0^*)$.
\end{definition} 
The meaning of this definition is that if $(x^*_0,\theta^*_0)$ are the true initial conditions of the agent, then a consistent posterior estimate is such that it collapses until all probability mass is on the true initial conditions. Statistical consistency of (\ref{eq:post_dens}) also needs an additional technical assumption:
\begin{assumption}
\label{ass:ass_7}
Let $(x^*_0,\theta^*_0)$ be the agent's true initial conditions.  The incentives $\pi_t$ are such that
\begin{equation}
\max_{\mathcal{E}(\delta)} \lim_{T\rightarrow\infty}\sum_{i=0}^{n_x} \log \frac{p_\nu(\tilde{x}_{t_i}-D\overline{x}_{t_i})}{p_\nu(\tilde{x}_{t_i}-Dx_{t_i})} + \sum_{j=0}^{n_u} \log \frac{ p_\omega(\tilde{y}_{\tau_j}-C\overline{u}_{\tau_j})}{p_\omega(\tilde{y}_{\tau_j}-Cu_{\tau_j})}= -\infty
\end{equation}
for any $\delta > 0$, almost surely, where $x_t,u_t$ are the states and decisions under initial conditions $(x^*_0,\theta^*_0)$, and $\overline{x}_t,\overline{u}_t$ are the states and decisions under initial conditions $(x_0,\theta_0)$.
\end{assumption}
This type of assumption is common in the adaptive control literature \citep{craig1987,astrom1995}, and is known as a \emph{sufficient excitation} or a \emph{sufficient richness} condition.  It is a mild condition because there are multiple ways of ensuring this condition holds \citep{bitmead1984,craig1987,astrom1995}.  One simple approach \citep{bitmead1984} is to compute an input $\pi_t$ and then add a small amount of random noise (whose value is known since it is generated by the coordinator) to the input before applying the input to the agent.

\begin{proposition}
If Assumptions \ref{ass:ass_1}--\ref{ass:ass_7} hold, then the estimated posterior distribution denoted by $\hat{p}(x_0,\theta_0|\{\tilde{x}_{t_i}\}_{i=0}^{n_x}, \{\tilde{y}_{\tau_i}\}_{i=0}^{n_u},\{\pi_i\}_{i=0}^T )$ and given in (\ref{eq:post_dens}) is consistent. \label{prop:consistent}
\end{proposition}
\proof{Proof: } Let $(x^*_0,\theta^*_0)$ be the agent's true initial conditions, and observe that
\begin{multline}
\log \hat{p}\big(x_0,\theta_0\big|\{\tilde{x}_{t_i}\}_{i=0}^{n_x}, \{\tilde{y}_{\tau_i}\}_{i=0}^{n_u},\{\pi_i\}_{i=0}^T\big) = \log \hat{p}\big(x^*_0,\theta^*_0\big|\{\tilde{x}_{t_i}\}_{i=0}^{n_x}, \{\tilde{y}_{\tau_i}\}_{i=0}^{n_u},\{\pi_i\}_{i=0}^T\big)
 + \\ \sum_{i=0}^{n_x} \log \frac{p_\nu(\tilde{x}_{t_i}-D\overline{x}_{t_i})}{p_\nu(\tilde{x}_{t_i}-Dx_{t_i})} + \sum_{j=0}^{n_u} \log \frac{ p_\omega(\tilde{y}_{\tau_j}-C\overline{u}_{\tau_j})}{p_\omega(\tilde{y}_{\tau_j}-Cu_{\tau_j})}+ \log \frac{p(x_0,\theta_0)}{p(x^*_0,\theta^*_0)} \label{eq:second_line},
\end{multline}
where $x_t,u_t$ are the states and decisions under initial conditions $(x^*_0,\theta^*_0)$, and $\overline{x}_t,\overline{u}_t$ are the states and decisions under initial conditions $(x_0,\theta_0)$.  But $\log \frac{p(x_0,\theta_0)}{p(x^*_0,\theta^*_0)}$ is a constant by assumption, and $\log \hat{p}\big(x^*_0,\theta^*_0\big|\{\tilde{x}_{t_i}\}_{i=0}^{n_x}, \{\tilde{y}_{\tau_i}\}_{i=0}^{n_u},\{\pi_i\}_{i=0}^T\big) \leq 0$ since $\hat{p}\big(x_0,\theta_0 \big|\{\tilde{x}_{t_i}\}_{i=0}^{n_x}, \{\tilde{y}_{\tau_i}\}_{i=0}^{n_u},\{\pi_i\}_{i=0}^T \big) \in[0,1]$ by construction.  So using Assumption \ref{ass:ass_6} gives $\max_{\mathcal{E}(\delta)}\log \hat{p}(x_0,\theta_0|\{\tilde{x}_{t_i}\}_{i=0}^{n_x}, \{\tilde{y}_{\tau_i}\}_{i=0}^{n_u},\{\pi_i\}_{i=0}^T) \rightarrow -\infty$ for any $\delta > 0$ almost surely.  Equivalently, $\max_{\mathcal{E}(\delta)}\hat{p}(x_0,\theta_0|\{\tilde{x}_{t_i}\}_{i=0}^{n_x}, \{\tilde{y}_{\tau_i}\}_{i=0}^{n_u},\{\pi_i\}_{i=0}^T) \rightarrow 0$ for any $\delta > 0$ almost surely. Thus for any $\delta > 0$ we have that
\begin{multline}
\hat{p}(\mathcal{E}(\delta)|\{\tilde{x}_{t_i}\}_{i=0}^{n_x}, \{\tilde{y}_{\tau_i}\}_{i=0}^{n_u},\{\pi_i\}_{i=0}^T) = \textstyle\int_{\mathcal{E}(\delta)}\hat{p}(x_0,\theta_0|\{\tilde{x}_{t_i}\}_{i=0}^{n_x}, \{\tilde{y}_{\tau_i}\}_{i=0}^{n_u},\{\pi_i\}_{i=0}^T)\times dx_0\times d\theta_0 \leq \\
\mathrm{volume}(\mathcal{X}\times\Theta)\cdot\textstyle\max_{\mathcal{E}(\delta)}\hat{p}(x_0,\theta_0|\{\tilde{x}_{t_i}\}_{i=0}^{n_x}, \{\tilde{y}_{\tau_i}\}_{i=0}^{n_u},\{\pi_i\}_{i=0}^T) \rightarrow 0 \label{eqn:deltaecal}
\end{multline}
almost surely.  This proves the result since (\ref{eqn:deltaecal}) holds almost surely for any $\delta > 0$.
\halmos 
\endproof
%This result implies statistical consistency of the MAP estimate.
\begin{corollary}
If Assumptions \ref{ass:ass_1}--\ref{ass:ass_7} hold, then $(\hat{x}_{0,T},\hat{\theta}_{0,T}) \stackrel{p}{\rightarrow} (x_0^*,\theta_0^*)$ as $T \rightarrow \infty$.\label{cor:map_consist}
%Let $(\hat{x}_{0,T},\hat{\theta}_{0,T}) = \argmax_{(x_0,\theta_0)} \psi_T(x_0,\theta_0)$ be the MAP estimate. If $p(x_0,\theta_0)>0$ for all $(x_0,\theta_0)\in \mathcal{X}\times\Theta$, and Assumptions \ref{ass:ass_1}--\ref{ass:ass_6} hold; then $(\hat{x}_{0,T},\hat{\theta}_{0,T}) \stackrel{p}{\rightarrow} (x_0^*,\theta_0^*)$ as $T \rightarrow \infty$.\label{cor:map_consist}
\end{corollary}
The above two results imply that future agent behavior can be reasonably predicted using the MAP parameters. Recall that calculating the MAP can be formulated as a MILP or MIQP, since the corresponding optimization problem is \eqref{eq:prof_lik_prob} with the constraints $x_0=\bar{x}_0,\theta_0=\bar{\theta}_0$ removed. 

%\section{Personalized Intervention Calculation}

\section{Optimizing Incentives}
 \label{sec:opt_sing_age}

The final step of our behavioral analytics framework involves using estimates of behavioral model parameters for each agent to optimize the design of costly incentives provided to the agents by the coordinator.   In Section \ref{sec:single}, we develop an algorithm for the single agent case.  In Section \ref{sec:multiple}, we use this single-agent algorithm as a sub-problem in the multi-agent case.  In both cases, we show that our algorithms are asymptotically optimal (as time continues and more data is collected) with respect to the coordinator's loss function when the agents behave according to the model constructed in Section \ref{sec:pred_mod}. % as part of the first step of our implementation of a behavioral analytics framework. 
The two algorithms we present in fact combine the second and third steps of our framework by first applying the parameter estimation algorithms (described in Section \ref{sec:pred_mod_est}) that comprise the second step, and then optimizing incentives.  The benefit of combining the second and third steps into a single algorithm is that this makes it easier to recompute the incentives as more data is collected over time from each agent. 

\subsection{Optimizing Incentives for a Single Agent} \label{sec:single}

Consider the problem of designing optimal incentives for a single agent at time $T$ by choosing $\{\pi_i\}_{i=T+1}^{T+n} \in \Pi^n$ to minimize a bounded loss function $\ell: \mathcal{X}^n\times \mathcal{U}^n \rightarrow \mathbb{R}$ of the agent's system states and decisions over the next $n$ time periods. We consider losses of a fairly general form:
\begin{assumption}
\label{ass:lossfun}
The loss function $\ell$ can be described by mixed integer linear constraints.
\end{assumption}
Since the coordinator only has noisy and incomplete observations of the agent's system states and decisions $\{\tilde{x}_{t_i}\}_{i=0}^{n_x}, \{\tilde{y}_{\tau_i}\}_{i=0}^{n_u}$, one design approach is to minimize the expected posterior loss
\begin{equation}
\label{eqn:eplsa}
\min\Big\{\mathbb{E}\big[\ell(\{x_t,u_t\}_{t=T+1}^{T+n})\big|\{\tilde{x}_{t_i}\}_{i=0}^{n_x}, \{\tilde{y}_{\tau_i}\}_{i=0}^{n_u}, \{\pi_i\}_{i=0}^T\big]\ \Big|\ \{\pi_i\}_{i=T+1}^{T+n} \in \Pi^n\Big\}.
\end{equation}
However, recalling our previous discussion, the agent's behavior is completely characterized by the initial conditions $(x_0,\theta_0)$, and so by the sufficiency and the smoothing theorem \citep{bickel2006}, there exists $\varphi: \mathcal{X}\times \Theta \times \Pi^n \mapsto \mathbb{R}$ such that the design problem can be exactly reformulated as
\begin{equation}
\min\Big\{\mathbb{E}\big[\varphi(x_0,\theta_0,\{\pi_i\}_{i=0}^{T+n})\big|\{\tilde{x}_{t_i}\}_{i=0}^{n_x}, \{\tilde{y}_{\tau_i}\}_{i=0}^{n_u}, \{\pi_i\}_{i=0}^T\big]\ \Big|\ \{\pi_i\}_{i=T+1}^{T+n} \in \Pi^n\Big\}.
\end{equation}
Calculating this expectation is difficult because the posterior distribution of $(x_{0},\theta_{0})$ does not generally have a closed form expression. In principle, discretization approaches from scenario generation \citep{kaut2003} could be used to approximate the design problem as 
\begin{equation}
\min\Big\{\frac{\sum_{i=1}^{M}\varphi({x}_{i,0},{\theta}_{i,0},\pi)\exp(\psi_T({x}_{i,0},{\theta}_{i,0};\pi))}{\sum_{i=1}^M  \exp(\psi_T({x}_{i,0},{\theta}_{i,0};\pi))}\ \Big|\ \{\pi_i\}_{i=T+1}^{T+n} \in \Pi^n\Big\}.
\label{eq:scen_gen_prob}
\end{equation}
where $(x_{i,0},\theta_{i,0})$ is an exhaustive enumeration of $\mathcal{X}\times\Theta$.  This approximation (\ref{eq:scen_gen_prob}) is still challenging to solve because the objective has a fractional, nonconvex form, and $\psi_T$ is defined as the value function of a MILP,  meaning that it does not have an easily computable closed form expression \citep{ralphs2014}. This means \eqref{eq:scen_gen_prob} is a BMIP with lower level problems that are MILP's. This is a complex class of optimization problems for which existing algorithms can only solve small problem instances \citep{moore1990,moore1992,denegre2009}.  

In this section, we develop a practical algorithm for optimizing incentives for a single agent.  We first summarize our algorithm, and show it only requires solving two MILP's.  Next we prove this algorithm can be interpreted as solving an approximation of solving either (\ref{eqn:eplsa}) or the optimal incentive design problem under perfect noiseless information.  More substantially, we also show that our algorithm provides a set of incentives that are asymptotically optimal as time advances.

\subsubsection{Two Stage Adaptive Algorithm (2SSA)}

Algorithm \ref{alg:2_ssa} summarizes our two stage adaptive approach (2SSA) for designing optimal incentives for a single agent.  The idea of the algorithm is to first compute a MAP estimate of the agent's initial conditions, use the MAP estimate as data for the first two arguments of $\varphi$, and then minimize $\varphi$.  In fact, we can solve this minimization problem without having to explicitly compute $\varphi$.  Because $\varphi$ is defined as the composition of the agent's dynamics with initial conditions $(x_0,\theta_0)$ and the coordinator's loss function $\ell$, it can be written as the value function of a feasibility problem:
\begin{equation}
\begin{aligned}
\varphi(\overline{x}_0,\overline{\theta}_0,\{\overline{\pi}_i\}_{i=0}^{T+n}) = \min_{x_t,u_t,\theta_t,\pi_t} & \ell(\{x_t,u_t\}_{t=T+1}^{T+n}) \\
 \text{s.t. }& u_t \in \argmax\big\{f(x_{t+1},u,\theta_t,\pi_t)\ \big|\ x_{t+1} = h(x_t,u),\ u\in\mathcal{U}\big\} & \forall t, \\
& x_{t+1} = h(x_t,u_t) & \forall t, \\
&\theta_{t+1} = g(x_t,u_t,\theta_t,\pi_t) &\forall t, \\
& x_t \in \mathcal{X}, \theta_t \in \Theta, \pi_t \in \Pi &\forall t, \\
& x_0 = \overline{x}_0, \; \theta_0 = \overline{\theta}_0, \pi_t = \overline{\pi}_t &\forall t.
\end{aligned}
\label{eq:phival_reform}
\end{equation}
More importantly, the problem of minimizing this $\varphi$ can be formulated as a MILP.
\begin{corollary} 
If Assumptions \ref{ass:ass_1}--\ref{ass:lossfun} hold, then $\varphi(x_0,\theta_0,\{\pi_i\}_{i=0}^{T+n})$ is lower semicontinuous in $x_0,\theta_0,\{\pi_i\}_{i=T+1}^{T+n}$, and the optimization problem given by $\min\{\varphi(x_0,\theta_0,\{\pi_i\}_{i=0}^{T+n})\ |\ \{\pi_i\}_{i=T+1}^{T+n} \in \Pi^n\}$ can be formulated as a MILP for all fixed values of $(x_0,\theta_0,\{\pi_i\}_{i=0}^T) \in \mathcal{X}\times \Theta\times\Pi^{T+1}$.
\label{cor:milp_reform}
\end{corollary} 
\proof{Proof: } For the first result, note Proposition \ref{prop:form_const} implies the feasible region of (\ref{eq:phival_reform}) can be expressed as mixed integer linear constraints with respect to $(x_t,u_t,\theta_t,\pi_t)$.  Thus $\varphi(\overline{x}_0,\overline{\theta}_0,\{\overline{\pi}_i\}_{i=0}^{T+n})$ is the value function of a MILP in which $x_0,\theta_0,\overline{\pi}_t$ belong to an affine term.  Standard results \citep{ralphs2014} imply the value function is lower semicontinuous with respect to $x_0,\theta_0,\{\pi_i\}_{i=T+1}^{T+n}$.

To show the second result, note that the problem of $\min\{\varphi(x_0,\theta_0,\{\pi_i\}_{i=0}^{T+n})\ |\ \{\pi_i\}_{i=T+1}^{T+n} \in \Pi^n\}$ is equivalent to (\ref{eq:phival_reform}) but with removal of the constraints $\pi_t = \overline{\pi}_t$ for $t = T+1,\ldots,T+n$.  And so the result follows by Proposition \ref{prop:form_const} and by recalling the assumptions on $\Pi$ and $\ell$.\halmos 
\endproof

\begin{algorithm}
\begin{algorithmic}[1]
\caption{Two Stage Single Agent Algorithm (2SSA)}\label{alg:2_ssa}
\Require $\{\tilde{x}_{t_i}\}_{i=0}^{n_x},\{\tilde{u}_{\tau_i}\}_{i=0}^{n_u}, \{\pi_i\}_{i=0}^T$
\State compute $(\hat{x}_{0,T},\hat{\theta}_{0,T}) \in \arg\max_{(x_0,\theta_0)} \psi_T(x_0,\theta_0)$
\State \Return $\pi_{2SSA}(T) \in \arg\min\{\varphi(\hat{x}_{0,T},\hat{\theta}_{0,T},\{\pi_i\}_{i=0}^{T+n})\ |\ \{\pi_i\}_{i=T+1}^{T+n} \in \Pi^n\}$
\end{algorithmic}
\end{algorithm}

\subsubsection{Asymptotic Optimality of 2SSA}
\label{prop:expect_conv}

The next result provides the underlying intuition of 2SSA.  In particular, we are approximating $\mathbb{E}[\varphi(x_0,\theta_0,\{\pi_i\}_{t=0}^{T+n})|\{\tilde{x}_{t_i}\}_{i=0}^{n_x}, \{\tilde{y}_{\tau_i}\}_{i=0}^{n_u}, \{\pi_i\}_{i=0}^T]$ using $\varphi(\hat{x}_{0,T},\hat{\theta}_{0,T},\{\pi_i\}_{i=0}^{T+n})$, and both these functions are converging to $\varphi(x_{0}^*,\theta_{0}^*,\{\pi_i\}_{t=0}^{T+n})$.
\begin{proposition}
Suppose that Assumptions \ref{ass:ass_1}--\ref{ass:lossfun} hold.  Then as $T\rightarrow\infty$ we have that: $\mathbb{E}[\varphi(x_0,\theta_0,\{\pi_i\}_{t=0}^{T+n})|\{\tilde{x}_{t_i}\}_{i=0}^{n_x}, \{\tilde{y}_{\tau_i}\}_{i=0}^{n_u}, \{\pi_i\}_{i=0}^T] \stackrel{p}{\rightarrow} \varphi(x_{0}^*,\theta_{0}^*,\{\pi_i\}_{t=0}^{T+n})$ for all fixed $\{\pi_i\}_{i=0}^{T+n}$; and $\varphi(\hat{x}_{0,T},\hat{\theta}_{0,T},\{\pi_i\}_{t=0}^{T+n}) \xrightarrow[\Pi^n]{\text{l-prob}} \varphi(x^*_{0},\theta^*_{0},\{\pi_i\}_{t=0}^{T+n})$. Here, $\Lambda_n \xrightarrow[\mathcal{X}]{\text{l-prob}} \Lambda$ means random function $\Lambda_n:\mathcal{X}\rightarrow\mathbb{R}$ is a lower semicontinuous approximation to function $\Lambda:\mathcal{X}\rightarrow\mathbb{R}$ \citep{vogel2003}.\label{prop:lscaf}
\end{proposition}
\proof{Proof: }
We begin by proving the first result.  By definition we have that
\begin{multline}
\mathbb{E}\big[\varphi(x_0,\theta_0,\{\pi_i\}_{i=0}^{T+n})\big|\{\tilde{x}_{t_i}\}_{i=0}^{n_x}, \{\tilde{y}_{\tau_i}\}_{i=0}^{n_u}, \{\pi_i\}_{i=0}^T\big] = \\\textstyle\int_{\mathcal{X}\times \Theta}\varphi(x_{0},\theta_{0},\{\pi_i\}_{i=0}^{T+n}) p(x_{0},\theta_{0}|\{\tilde{x}_{t_i}\}_{i=0}^{n_x}, \{\tilde{y}_{\tau_i}\}_{i=0}^{n_u}, \pi)\times dx_0\times d\theta_0. \label{eq:line1_prop4}
\end{multline} 
Also, Proposition \ref{prop:consistent} implies the posterior $p(x_{0},\theta_{0}|\{\tilde{x}_{t_i}\}_{i=0}^{n_x}, \{\tilde{y}_{\tau_i}\}_{i=0}^{n_u}, \{\pi_i\}_{i=0}^{T})$ is consistent, and thus becomes degenerate at $(x_{0}^*,\theta_{0}^*)$ in the limit. Hence the Dominated Convergence Theorem gives
\begin{equation}
\eqref{eq:line1_prop4}\stackrel{p}{\rightarrow} \textstyle\int_{\mathcal{X}\times \Theta}\varphi(x_{0},\theta_{0},\pi)\times \delta(x_0 -x_0^*)\times\delta(\theta_0 -\theta_0^*)\times  dx_0\times d\theta_0 = \varphi(x_{0}^*,\theta_{0}^*,\pi),
\end{equation}
where in the equation above $\delta(\cdot)$ is the Dirac delta function.

For the second result, recall that Corollary \ref{cor:map_consist} implies $(\hat{x}_{0,T},\hat{\theta}_{0,T})\stackrel{p}{\rightarrow}(x_0^*,\theta_0^*)$.  And Corollary \ref{cor:milp_reform} gives that $\varphi(x_0,\theta_0,\{\pi_i\}_{i=0}^{T+n})$ is lower semicontinuous in $x_0,\theta_0,\{\pi_i\}_{i=T+1}^{T+n}$.  The result then follows by direct application of Proposition 2.1.ii of \citep{vogel2003b}.\halmos
\endproof
If the coordinator had perfect knowledge of the agent's true initial conditions $(x_0^*,\theta_0^*)$, then the optimal incentives are $\arg\min\{\varphi(x_0^*,\theta_0^*,\{\pi_i\}_{i=0}^{T+n})\ |\ \{\pi_i\}_{i=T+1}^{T+n} \in \Pi^n\}$.  But since we do not know the initial conditions, the above result shows that both (\ref{eqn:eplsa}) and $\arg\min\{\varphi(\hat{x}_{0,T},\hat{\theta}_{0,T},\{\pi_i\}_{i=0}^{T+n})\ |\ \{\pi_i\}_{i=T+1}^{T+n} \in \Pi^n\}$ are reasonable approximations.  In fact, we can show a stronger result for the solution generated by 2SSA.
\begin{theorem}
Note that $\arg\min\{\varphi(x^*_{0},\theta^*_{0},\{\pi_i\}_{i=0}^{T+n}) | \{\pi_i\}_{i=T+1}^{T+n}\in\Pi^n\}$ is the set of optimal solutions under the agent's true initial conditions $(x^*_0,\theta^*_0)$. If Assumptions \ref{ass:ass_1}--\ref{ass:lossfun} hold, then we have that
\begin{equation}
\mathrm{dist}\Big(\pi_{2SSA}(T), \arg\min\{\varphi(x_0^*,\theta_0^*,\{\pi_i\}_{i=0}^{T+n})\ |\ \{\pi_i\}_{i=T+1}^{T+n} \in \Pi^n\}\Big) \stackrel{p}{\rightarrow} 0
\end{equation}
as $T \rightarrow \infty$, for any $\pi_{2SSA}(T)$ returned by 2SSA.  Note $\mathrm{dist}(x,B) = \inf_{y\in B}\|x-y\|$.
\label{thm:alg_opti}
\end{theorem}
\proof{Proof: }
The result follows by combining the second part of our Proposition \ref{prop:lscaf} with Theorem 4.3 from \citep{vogel2003}. \halmos
\endproof
This result suggests that any solution returned by 2SSA is asymptotically included within the set of optimal incentives computed for the agent's true initial conditions.  Restated, the result says 2SSA provides a set of incentives that are asymptotically optimal.  This is a non-obvious result because in general pointwise-convergence of a sequence of stochastic optimization problems is not sufficient to ensure convergence of the minimizers of the sequence of optimization problems to the minimizer of the limiting optimization problem.  \cite{rockafellar2009variational} provide an example that demonstrates this possible lack of convergence of minimizers.

\subsection{Policy Calculation With Multiple Agents} \label{sec:multiple}
%\label{sec:mult_aget}
We next study the general setting where the coordinator designs incentives for a large group of agents. We let $\mathcal{A}$ be the set of agents, and the quantities corresponding to a specific agent $a\in\mathcal{A}$ are denoted using subscript $a$.  Now suppose that at time $T$ the coordinator measures $\{\tilde{x}^a_{t_i}\}_{i=0}^{n^a_x}, \{\tilde{y}^a_{\tau_i}\}_{i=0}^{n^a_u}$ for all agents $a \in \mathcal{A}$. One approach to designing incentives is by solving:
\begin{multline} 
\min \Big\{\mathbb{E}\big[\Phi(x^a_{0},\theta^a_{0},\{\pi_i^a\}_{i=0}^{T+n} \text{ for } a \in \mathcal{A})|\{\tilde{x}^a_{t_i}\}_{i=0}^{T^a_x}, \{\tilde{y}^a_{\tau_i}\}_{i=0}^{T^a_u}, \{\pi_i^a\}_{i=0}^{T} \text{ for } a \in \mathcal{A}\big]\ \Big|\\\big\{ \{\pi_i^a\}_{i=T+1}^{T+n}\text{ for } a \in \mathcal{A}\big\} \in \Omega \Big\}\label{eq:dec_mult_prob}
\end{multline}
Here, $\Phi:\mathcal{X}^{\#\mathcal{A}}\times\theta^{\#\mathcal{A}}\times\Omega\rightarrow \mathbb{R}$ is a joint loss function that depends on the behavior all agents.  For the settings we are interested in, this loss function has a separable structure.
%, and $\pi$ is a matrix such that each of its columns $\pi^a$ corresponds to the strategy used for agent $a \in \mathcal{A}$. We make the following additional assumptions to analyze this problem:
\begin{assumption} \label{ass:seperab}
Loss $\Phi$ is additively $\Phi(x^a_{0},\theta^a_{0},\{\pi_i^a\}_{i=0}^{T+n}\ \mathrm{ for }\ a \in \mathcal{A}) = \sum_{a \in \mathcal{A}} \varphi^a(x^a_{0},\theta^a_{0},\{\pi_i^a\}_{i=0}^{T+n})$ or multiplicatively separable $\Phi(x^a_{0},\theta^a_{0},\{\pi_i^a\}_{i=0}^{T+n}\ \mathrm{ for }\ a \in \mathcal{A}) = \prod_{a \in \mathcal{A}} \varphi^a(x^a_{0},\theta^a_{0},\{\pi_i^a\}_{i=0}^{T+n})$.%, where $\varphi_a$ are defined as in \eqref{eq:phival_reform} but with loss function $\ell_a$.% with respect to each agent's parameters and policy.
\end{assumption} 
Without loss of generality, we assume $\Phi$ is additively separable since we can obtain similar results for the case of multiplicative separability by taking the logarithm of $\Phi$.  We also make an assumption that states $\Omega$ is decomposable in a simple way.

%We also make an assumption about the structure of $\Pi^{|\mathcal{A}|}$.
%Where $\varphi_a$ in Assumption \ref{ass:seperab} are defined as functional compositions in the same sense as Section \ref{sec:opt_sing_age}. Without loss of generality, we will consider that $\Phi$ is additively separable for our analysis, since we can obtain similar results for the multiplicative separability case by taking the logarithm of $\Phi$. Next we consider assumptions on the structure of $\Pi^{|\mathcal{A}|}$.
\begin{assumption}
\label{ass:svarep}
There exist a finite set $V = \{v_1,v_2,\ldots\}$ with vector-valued, sets $S_{v} \subseteq \Pi^n$ for $v\in V$, and a vector-valued constant $\beta$ such that
\begin{multline}
\Omega = \Big\{\{\pi_i^a\}_{i=T+1}^{T+n}\ \mathrm{ for }\ a \in \mathcal{A} : y_v^a\in\{0,1\}, \textstyle\sum_{v\in V}y_v^a = 1\ \mathrm{ for }\ a \in\mathcal{A}, \textstyle\sum_{a\in\mathcal{A}}\sum_{v\in V} v\cdot y_{v}^a \leq \beta\\ \{\pi_i^a\}_{i=T+1}^{T+n} \in S_{v}\ \mathrm{if}\ y_v^a = 1\Big\}.
\end{multline}
Moreover, the sets $S_v$ are compact sets that are representable by a finite number of mixed integer linear constraints.
\end{assumption}
The underlying idea of this assumption is that the set $V$ describes a vector of discrete elements that can be used as incentives, and $\beta$ is the vector-valued budget on the discrete incentives.  When the discrete incentives are fixed at $v$, the set $S_{v}$ keeps the discrete incentives fixed and describes the feasible set of continuous incentives.

Even with these assumptions on separability and decomposibility, solving \eqref{eq:dec_mult_prob} is difficult because it is a BMIP with $\# \mathcal{A}$ MILP's in the lower level.  Thus, we develop an adaptive algorithm (based on the 2SSA algorithm) for optimizing incentives for multiple agents.  We first summarize our algorithm, and demonstrate that it only requires solving a small number of computable MILP's.  Next we prove this algorithm provides a set of incentives that are asymptotically optimal as time advances.

% In this section, we develop a practical algorithm for optimizing incentives for a single agent.  We first summarize our algorithm, and show it only requires solving two MILP's.  Next we prove this algorithm can be interpreted as solving an approximation of solving either (\ref{eqn:eplsa}) or the optimal incentive design problem under perfect noiseless information.  More substantially, we also show that our algorithm provides a set of incentives that are asymptotically optimal as time advances.

\subsubsection{Adaptive Algorithm for Multiple Agents}

\begin{algorithm}
\caption{Adaptive Behavioral Multi-Agent Algorithm (ABMA)}\label{alg:mult_age_alg}
\begin{algorithmic}[1]
\Require $\{\tilde{x}_{t_i}^a\}_{i=0}^{n_x^a},\{\tilde{u}_{\tau_i}^a\}_{i=0}^{n_u^a}, \{\pi_i^a\}_{i=0}^T$ for $a\in\mathcal{A}$
\ForAll{$a \in \mathcal{A}$}
\State compute $(\hat{x}^{a}_{0,T},\hat{\theta}^{a}_{0,T}) = \argmax_{(x_0,\theta_0)} \psi_T(x_0,\theta_0)$
\ForAll{$v\in V$}
\State set $\pi^a_v \in \arg\min\{\varphi^a(\hat{x}_{0,T},\hat{\theta}_{0,T},\{\pi_i\}_{i=0}^{T+n})\ |\ \{\pi_i\}_{i=T+1}^{T+n} \in S_v\}$ 
\State set $\phi^a_v =  \varphi^a(\hat{x}^{a}_{0,T},\hat{\theta}^{a}_{0,T},\pi^a_v)$
\EndFor
\EndFor
\State compute:
\begin{equation*}
\begin{aligned}
y \in \textstyle\argmin\ &\textstyle\sum_{a\in \mathcal{A}}\sum_{v\in V} \phi^a_v\cdot y^a_v\\
\text{s.t. }&\textstyle\sum_{a\in\mathcal{A}}\sum_{v\in V} v\cdot y_{v}^a \leq \beta\\
&\textstyle\sum_{v\in V}y_v^a = 1\ \mathrm{ for }\ a \in\mathcal{A}\\
&y_v^a\in\{0,1\}
\end{aligned}
\end{equation*}
\ForAll{$a \in \mathcal{A}$ and $v \in V$}
\State{set $\pi^a_{ABMA}(T) = \pi^{a}_{v}$} \textbf{if} $y^a_v = 1$
\EndFor
\State \Return $\pi^a_{ABMA}(T)$ for $a\in\mathcal{A}$
%\EndProcedure
\end{algorithmic}
\end{algorithm}

We design incentives for multiple agents with the Adaptive Behavioral Multi-Agent Algorithm (ABMA) presented in Algorithm \ref{alg:mult_age_alg}. The main idea behind this method is to use the assumptions on $\Phi$ and $\Omega$ to decompose the initial problem into $\#\mathcal{A}$ sub-problems that solve a single agent problem, and a single master problem that combines these solutions into a global optimum across all agents. Because of the assumptions on $\Omega$, each sub-problem can be further decomposed into $\#V$ sub-problems.  For each sub-problem, we use the 2SSA algorithm to solve the $\#\mathcal{A}\cdot\#V$ sub-problems; however, we do not explicitly call the 2SSA algorithm because it is more efficient to solve the MAP estimator once and then solve the incentive design problem for each single agent.  Our first result concerns the computability of this algorithm.

\begin{proposition}
If Assumptions \ref{ass:ass_1}--\ref{ass:svarep} hold, then the main computational steps of the ABMA algorithm involve solving a total of $\#\mathcal{A}\cdot(\#V + 1)$ MILP's and 1 ILP.
\end{proposition}

\proof{Proof: } 
Step 2 of ABMA is a MAP estimate, which can be computed by solving a single MILP by Corollary \ref{cor:mapmilp}.  A similar argument used to prove Corollary \ref{cor:milp_reform} shows that Steps 4 and 5 can be computed by solving a single MILP.  Step 8 can be seen to be an ILP by construction.  The remaining steps of ABMA are assignment steps and do not require solving any optimization problems.
\halmos 
\endproof
This means the ABMA algorithm performs incentive design for the multi-agent case by solving $\#\mathcal{A}\cdot(\#V+1) + 1$ MILP's, which is significantly less challenging than solving a BMIP with $\#\mathcal{A}$ MILP's in the lower level as would be required to solve \eqref{eq:dec_mult_prob}.		

The ABMA algorithm also has an alternative interpretation, and to better understand this consider the following feasibility problem:
\begin{multline}
\label{eqn:abmaundec}
\Phi(\overline{x}^a_{0},\overline{\theta}^a_{0},\{\overline{\pi}_i^a\}_{i=0}^{T+n} \text{ for } a \in \mathcal{A}) = \\
\begin{aligned}
\min_{x_t^a,u_t^a,\theta_t^a,\pi_t^a} & \Phi(x^a_{0},\theta^a_{0},\{\pi_i^a\}_{i=0}^{T+n} \text{ for } a \in \mathcal{A}) \\
 \text{s.t. }& u_t^a \in \argmax\big\{f(x_{t+1}^a,u,\theta_t^a,\pi_t^a)\ \big|\ x_{t+1}^a = h(x_t^a,u),\ u\in\mathcal{U}\big\} & \forall a,t, \\
& x_{t+1}^a = h(x_t^a,u_t^a) & \forall a,t, \\
&\theta_{t+1}^a = g(x_t^a,u_t^a,\theta_t^a,\pi_t^a) &\forall a,t, \\
& x_t^a \in \mathcal{X}^a, \theta_t^a \in \Theta, \pi_t^a \in \Pi &\forall a,t, \\
& x_0^a = \overline{x}_0^a, \theta_0^a = \overline{\theta}^a_0, \{\pi_t^a\}_{t=0}^{T+n} = \{\overline{\pi}_t^a\}_{t=0}^{T+n} &\forall a,t,\\
&\{ \{\pi_t^a\}_{t=T+1}^{T+n}\ \mathrm{ for }\ a \in \mathcal{A}\} \in \Omega.
\end{aligned}
\end{multline}
Our first result concerns regularity properties of the above written feasibility problem.

\begin{proposition}
\label{prop:lscphi}
If Assumptions \ref{ass:ass_1}--\ref{ass:svarep} hold, then $\Phi(x^a_{0},\theta^a_{0},\{\pi_i^a\}_{i=0}^{T+n}\ \mathrm{ for }\ a \in \mathcal{A})$ is lower semicontinuous in its arguments, and $\min\{\Phi(x^a_{0},\theta^a_{0},\{\pi_i^a\}_{i=0}^{T+n}\ \mathrm{ for }\ a \in \mathcal{A})\ |\ \{ \{\pi_t^a\}_{t=T+1}^{T+n}\ \mathrm{ for }\ a \in \mathcal{A}\} \in \Omega\}$ can be formulated as a MILP for all fixed values of $(x_0^a,\theta_0^a,\{\pi_i^a\}_{i=0}^T) \in \mathcal{X}\times \Theta\times\Pi^{T+1}$ for $a\in\mathcal{A}$.
\end{proposition}

\proof{Proof: } 
Using the assumptions on separability of the joint loss function $\Phi$ (Assumption \ref{ass:seperab}) and decomposibility on the incentive set $\Omega$ (Assumption \ref{ass:svarep}), we have that (\ref{eqn:abmaundec}) can be reformulated as
\begin{equation}
\label{eqn:refundec}
\begin{aligned}
\min_{y_v^a}\ & \textstyle\sum_{a\in\mathcal{A}}\sum_{v\in V} \phi^a_v\cdot y^a_v\\
\text{s.t. } & \min\{\varphi^a({x}_{0,T},{\theta}_{0,T},\{\pi_i\}_{i=0}^{T+n})\ |\ \{\pi_i\}_{i=T+1}^{T+n} \in S_v\} \leq \phi^a_v\\
& x_0^a = \overline{x}_0^a, \theta_0^a = \overline{\theta}^a_0, \{\pi_t^a\}_{t=0}^{T+n} = \{\overline{\pi}_t^a\}_{t=0}^{T+n} \text{ for all } a,t\\
&\textstyle\sum_{a\in\mathcal{A}}\sum_{v\in V} v\cdot y_{v}^a \leq \beta\\
&\textstyle\sum_{v\in V}y_v^a = 1\ \mathrm{ for }\ a \in\mathcal{A}\\
&y_v^a\in\{0,1\}
\end{aligned}
\end{equation}
Since $\phi^a_v\cdot y^a_v$ is the product of a continuous and binary decision variable, standard integer programming reformulation techniques allow us to reformulate the above as
\begin{equation}
\begin{aligned}
\min_{y_v^a}\ & \textstyle\sum_{a\in\mathcal{A}}\sum_{v\in V} z^a_v\\
\text{s.t. } & \varphi^a({x}_{0,T},{\theta}_{0,T},\{\pi_i\}_{i=0}^{T+n}) \leq \phi^a_v\\
&\{\pi_i\}_{i=T+1}^{T+n} \in S_v \\
& x_0^a = \overline{x}_0^a, \theta_0^a = \overline{\theta}^a_0, \{\pi_t^a\}_{t=0}^{T+n} = \{\overline{\pi}_t^a\}_{t=0}^{T+n} \text{ for all } a,t\\
&\textstyle\sum_{a\in\mathcal{A}}\sum_{v\in V} v\cdot y_{v}^a \leq \beta\\
&\textstyle\sum_{v\in V}y_v^a = 1\ \mathrm{ for }\ a \in\mathcal{A}\\
&y_v^a\in\{0,1\}\\
&z^a_v \geq \phi^a_v - M\cdot(1-y^a_v)\\
&z^a_v \leq \phi^a_v + M\cdot(1-y^a_v)\\
&z^a_v \leq M\cdot y^a_v\\
&z^a_v \geq -M\cdot y^a_v\\
\end{aligned}
\end{equation}
where $M > 0$ is a large-enough constant.  Such a finite $M$ exists because $\mathcal{X},\mathcal{U},\Pi,\Theta$ are compact, and because $\ell^a$ is representable by a finite number of mixed integer linear constraints.  Since Corollary \ref{cor:milp_reform} (and its proof) implies we can represent $\varphi^a({x}_{0,T},{\theta}_{0,T},\{\pi_i\}_{i=0}^{T+n}) \leq \phi^a_v$ and $\{\pi_i\}_{i=T+1}^{T+n} \in S_v$ by mixed integer linear constraints, this means we can reformulate (\ref{eqn:abmaundec}) as a MILP with linear constraints that are affine in $(\overline{x}^a_{0},\overline{\theta}^a_{0},\{\overline{\pi}_i^a\}_{i=0}^{T+n} \text{ for } a \in \mathcal{A})$.  And so standard results \citep{ralphs2014} imply its value function is lower semicontinuous with respect to these variables, which is our first result.  The second result follows by noting $\min\{\Phi(x^a_{0},\theta^a_{0},\{\pi_i^a\}_{i=0}^{T+n}\ \mathrm{ for }\ a \in \mathcal{A})\ |\ \{ \{\pi_t^a\}_{t=T+1}^{T+n}\ \mathrm{ for }\ a \in \mathcal{A}\} \in \Omega\}$ is equivalent to (\ref{eqn:abmaundec}) but with removal of the constraints $\pi_t^a = \overline{\pi}_t^a$ for $t = T+1,\ldots,T+n$.\halmos 
\endproof

The optimization problem (\ref{eqn:abmaundec}) and the above result provide an alternative interpretation of the ABMA algorithm, which is formalized by the next corollary.

\begin{corollary}
\label{corr:altint}
If Assumptions \ref{ass:ass_1}--\ref{ass:svarep} hold, then the solution of $\min\{\Phi(x^a_{0},\theta^a_{0},\{\pi_i^a\}_{i=0}^{T+n} \text{ for } a \in \mathcal{A})\ |\ \{ \{\pi_t^a\}_{t=T+1}^{T+n}\text{ for } a \in \mathcal{A}\} \in \Omega\}$, is given by the ABMA algorithm but with Step 2 replaced with the step: $\mathrm{set}\ (\hat{x}^{a}_{0,T},\hat{\theta}^{a}_{0,T}) = (x_0^a,\theta_0^a)$.
\end{corollary}

\proof{Proof: } 
This is straightforward from the reformulation shown in (\ref{eqn:refundec}).\halmos \endproof
Thus, though (\ref{eqn:abmaundec}) is a large MILP, the assumptions we have made allow us to decompose the solution of this problem into a series of substantially smaller MILP's.

\subsubsection{Asymptotic Optimality of ABMA}

The optimization problem in (\ref{eqn:abmaundec}) is a useful construction because it can also be used to compute the optimal set of incentives.  If each agent's true initial conditions $(x_0^{*,a},\theta_0^{*,a})$ were known, then an optimal solution belongs to $\arg\min\{\Phi(x^{*,a}_{0},\theta^{*,a}_{0},\{\pi_i^a\}_{i=0}^{T+n} \text{ for } a \in \mathcal{A})\ |\ \{ \{\pi_t^a\}_{t=T+1}^{T+n}\text{ for } a \in \mathcal{A}\} \in \Omega\}$.  More importantly, we have the following relationship to the solutions of the ABMA algorithm:

\begin{theorem}
Note that $\arg\min\{\Phi(x^{*,a}_{0},\theta^{*,a}_{0},\{\pi_i^{*,a}\}_{i=0}^{T+n} \text{ for } a \in \mathcal{A})\ |\ \{ \{\pi_t^a\}_{t=T+1}^{T+n}\text{ for } a \in \mathcal{A}\} \in \Omega\}$ is the set of optimal solutions under the agents' true initial conditions $(x^{*,a}_0,\theta^{*,a}_0)$. If Assumptions \ref{ass:ass_1}--\ref{ass:lossfun} hold, then we have that
\begin{multline}
\mathrm{dist}\Big(\{\pi_{ABMA}^a(T)\ \mathrm{for}\ a\in\mathcal{A}\}, \\\arg\min\{\Phi(x^{*,a}_{0},\theta^{*,a}_{0},\{\pi_i^a\}_{i=0}^{T+n} \text{ for } a \in \mathcal{A})\ |\ \{ \{\pi_t^a\}_{t=T+1}^{T+n}\text{ for } a \in \mathcal{A}\} \in \Omega\}\Big) \stackrel{p}{\rightarrow} 0
\end{multline}
as $T \rightarrow \infty$, for any $\pi_{ABMA}^a(T)$ returned by ABMA.  Recall that $\mathrm{dist}(x,B) = \inf_{y\in B}\|x-y\|$.
\label{thm:alg_multiopti}
\end{theorem}

\proof{Proof: }
Corollary \ref{cor:map_consist} implies $(\hat{x}_{0,T}^a,\hat{\theta}_{0,T}^a)\stackrel{p}{\rightarrow}(x_0^{*,a},\theta_0^{*,a})$, and Corollary \ref{cor:milp_reform} states $\Phi$ is lower semicontinuous in its arguments. This means $\Phi(\hat{x}^{a}_{0},\hat{\theta}^{a}_{0},\{\pi_i^{a}\}_{i=0}^{T+n} \text{ for } a \in \mathcal{A})$ is a lower semicontinuous approximation to $\Phi(x^{*,a}_{0},\theta^{*,a}_{0},\{\pi_i^{a}\}_{i=0}^{T+n} \text{ for } a \in \mathcal{A})$ by Proposition 2.1.ii of \citep{vogel2003b}.  But Corollary \ref{corr:altint} shows that
\begin{equation}
\{\pi_{ABMA}^a(T)\ \mathrm{for}\ a\in\mathcal{A}\} \in \arg\min\{\Phi(\hat{x}^{a}_{0},\hat{\theta}^{a}_{0},\{\pi_i^a\}_{i=0}^{T+n} \text{ for } a \in \mathcal{A})\ |\ \{ \{\pi_t^a\}_{t=T+1}^{T+n}\text{ for } a \in \mathcal{A}\} \in \Omega\}.
\end{equation}
This means that the result follows by applying Theorem 4.3 from \citep{vogel2003}. \halmos
\endproof
Thus any solution returned by ABMA is asymptotically included within the set of optimal incentives computed for the agents' true initial conditions.  Restated, the above result says ABMA provides a set of incentives that are asymptotically optimal.  This is a non-obvious result because in general pointwise-convergence of a sequence of stochastic optimization problems is not sufficient to ensure convergence of the minimizers of the sequence of optimization problems to the minimizer of the limiting optimization problem.  \cite{rockafellar2009variational} provide an example that demonstrates this possible lack of convergence of minimizers.

\section{Computational Experiments: Weight Loss Program Design}\label{sec:experiment}
We have completed computational experiments applying the tools and techniques developed in this paper that form a specific implementation of a behavioral analytics framework.  We compare several approaches, including ours, for designing incentives for multiple myopic agents %our approachdifferent approaches (including our ABMA algorithm) 
to the problem described in Section \ref{sec:health_app} of designing personalized behavioral incentives for a clinically-supervised weight loss program.  The first step of our behavioral analytics approach is to construct a behavioral model of individuals in weight loss programs.  We describe the data source used for the simulations, and then summarize our behavioral model \citep{aswani2016} for individuals participating in such loss programs.  To demonstrate the second and third steps of our framework, we simulate a setting in which behavioral incentives chosen using our ABMA algorithm are evaluated against behavioral incentives computed by (intuitively-designed) adaptive heuristics.  Both our implementation of a behavioral analytics framework and the heuristic provide adaptation by recomputing the incentives at regular intervals as more data is collected from each individual.  Our metric for comparison is the number of individuals who achieve clinically significant weight loss (i.e., a 5\% reduction in body weight) at the end of the program.  We also compare the percentage of weight loss for individuals who do not achieve clinically significant weight loss in order to better understand how clinical visits are allocated by the different methods.  We conclude by performing a sensitivity analysis of design choices for the second and third steps of our behavioral analytics framework. 

\subsection{mDPP Program Trial Data Source}
Our computational experiments used data from the mDPP trial \citep{fukuoka2014b}. This was a randomized control trial (RCT) that was conducted to evaluate the efficacy of a 5 month mobile phone based weight loss program among overweight and obese adults at risk for developing type 2 diabetes. This program was adapted from the Diabetes Prevention Program (DPP) \citep{dpp2002,dpp2009}, but the number of in person clinical visits was reduced from 16 to 6 per person, and group exercise sessions were replaced with a home based exercise program to reduce costs. Sixty one overweight adults were randomized into an active control group which only received an accelerometer (n=31) or a treatment group which receive the mDPP mobile app plus the accelerometer and clinical office visits (n=30). Changes in primary and secondary outcomes for the trial were promising. The treatment group lost an average of 6.2 $\pm$ 5.9 kg (-6.8\% $\pm$ 5.7\%) between baseline and the 5 month follow up while the control group gained 0.3 $\pm$ 3.0 kg (0.3\% $\pm$ 5.7 \%) (p $<$ 0.001). The treatment group's steps per day increased by 2551 $\pm$ 4712 compared to the control group's decrease of 734 $\pm$ 3308 steps per day (p $<$ 0.001). Additional details on demographics and other treatment parameters are available in \citep{fukuoka2014b}.  The data available from the mDPP trial includes step data (from accelerometer measurements), body weight data (which was measured at least twice a week every week and recorded in the mobile app by individuals in the treatment group, as well as measured three times in a clinical setting at  baseline, 3 month, and 5 month), and demographic data (i.e., age, gender, and height of each individual). We note that this data matches the assumptions in Section \ref{sec:pred_mod}.

\subsection{Summary of Behavioral Model}

We construct a behavioral model for each individual participating in the weight loss program.  Using the terminology and notation of Section \ref{sec:pred_mod}, the system state of each individual $x_t$ is their body weight on day $t$, and their decisions $u_t = (f_t,s_t)$ on day $t$ are how many calories they consume $f_t$ and how many steps they walk $s_t$.  The behavioral incentives $\pi_t =(g_t,d_t)$ provided to an individual on day $t$ consist of (numeric) step goals $g_t$ and an indicator $d_t$ equal to one if a clinical visit was scheduled for that day.  The motivational state (or type) is $\theta_t := (u_b,f_{b,t},F_{b,t},p_t,\mu,\delta,\beta)$. The state $u_b$ captures the individual's baseline preference for number of steps taken each day, while $f_{b,t},F_{b,t}$ capture the individual's short term and long term caloric intake preference, respectively. The variable $p_t$ captures the disutility an individual experiences from not meeting a step goal. The last set of motivational states describe the individual's response to behavioral incentives. The states $\beta,\delta$ describe the amount of change in the individual's caloric consumption and physical activity preferences, respectively, after undergoing a single clinical visit. The state $\mu$ describes the self efficacy effect \citep{bandura1998,conner1996} from meeting exercise goals. 
 
 The utility function of an individual on day $t$ is given by $f(w_{t+1},f_t,u_t;\theta_t,\pi_t) = -w_{t+1}^2 -r(u - u_b)^2 -r(f - f_{b,t})^2 + p_t(u_t-g_t)^-$.  In our past modeling work \citep{aswani2016}, we found that the predictions of this behavioral model were relatively insensitive to the value of $r$.  And our numerical experiments in \citep{aswani2016} found that fixing the value of $r$ to be the same for each individual provided accurate predictions.  And so we assume that $r$ is a fixed and known constant in our numerical experiments here. The temporal dynamics for an individual's system and motivational states are 
 the dynamics of an individual's type by:
 \begin{align}
 w_{t+1} &= a\cdot w_t + b\cdot u_t + f_t + k \label{eq:wet_dynam}\\
 F_{b,t+1} &= (1-\alpha)\cdot F_{b,t} + \alpha\cdot f_{b,t} \label{eq:bas_food_dyn}\\
 f_{b,t+1} &= \gamma\cdot (f_{b,t} - F_{b,t}) + F_{b,t} - \beta\cdot d_n \label{eq:food_dyn}\\
 p_{t+1} &= \gamma\cdot p_t + \delta\cdot d_t + \mu\cdot \mathds{1}(u_t \geq g_t). \label{eq:step_dyn}
 \end{align} 
Equation \eqref{eq:wet_dynam} is a ``calories in minus calories out'' description of weight change, and a standard physiological formula \citep{mifflin1990} is used to compute the values of $a,b,k$ based on the demographics of the individual. Equation \eqref{eq:bas_food_dyn} models the long term caloric intake preference as an exponential moving average of the short term caloric intake preferences.  We found that the predictions for different individuals were relatively insensitive to the value of $\alpha$, and so in our numerical experiments we assume $\alpha$ is known and fixed to a value satisfying $\alpha < 1$. In \eqref{eq:food_dyn}, we model the dynamics of baseline food consumption as always tending towards their initial value unless perturbed by a clinical visit. In \eqref{eq:step_dyn} we model the tendency for meeting the step goal as tending towards zero unless there is a clinical visit or the individual has met the previous exercise goal, which increases their self efficacy and makes the individual more likely to meet their step goal in the future. In both \eqref{eq:food_dyn} and \eqref{eq:step_dyn}, $\gamma <0$ is assumed to be a known decay factor since we found that predictions were relatively insensitive to the value of $\gamma$ \citep{aswani2016}.  Note that these temporal dynamics and utility functions satisfy the assumptions in Section \ref{sec:pred_mod}.

For the MLE and MAP calculations, we assumed that step and weight data were measured with zero-mean noise distributed according to a Laplace distribution with known variance.  We found that predictions of the behavioral model estimated when assuming the noise had a Gaussian distribution were of the same quality as those estimated when assuming Laplace noise, and so we assume Laplace noise so that the MAP and MLE problems can be formulated as MILP's (as shown in Section \ref{sec:pred_mod}). Furthermore, the prior distribution used for the MAP calculation for each individual was a histogram of the MLE estimates of all the other individual's parameters. Note that this form for a prior distribution can be expressed using integer linear constraints \citep{aswani2016}. The complete MILP formulations for MAP and MLE are provided in the appendix.  
 
\subsection{Weight Loss Program Design}

Since the majority of implementation costs for weight loss programs are due to clinical visits, the clinician's design problem is to maximize the expected number of individuals who reduce their weight by a clinically significant amount (i.e., 5\% reduction in body weight).  The clinician is able to personalize the step goals for each individual, and can change the number and timing of clinical visits for each individual.  However, there is a budget constraint on the total number of visits that can be scheduled across all individuals.  This constraint captures the costliness of clinical visits.

\begin{figure}
\centering
\includegraphics[scale=0.8]{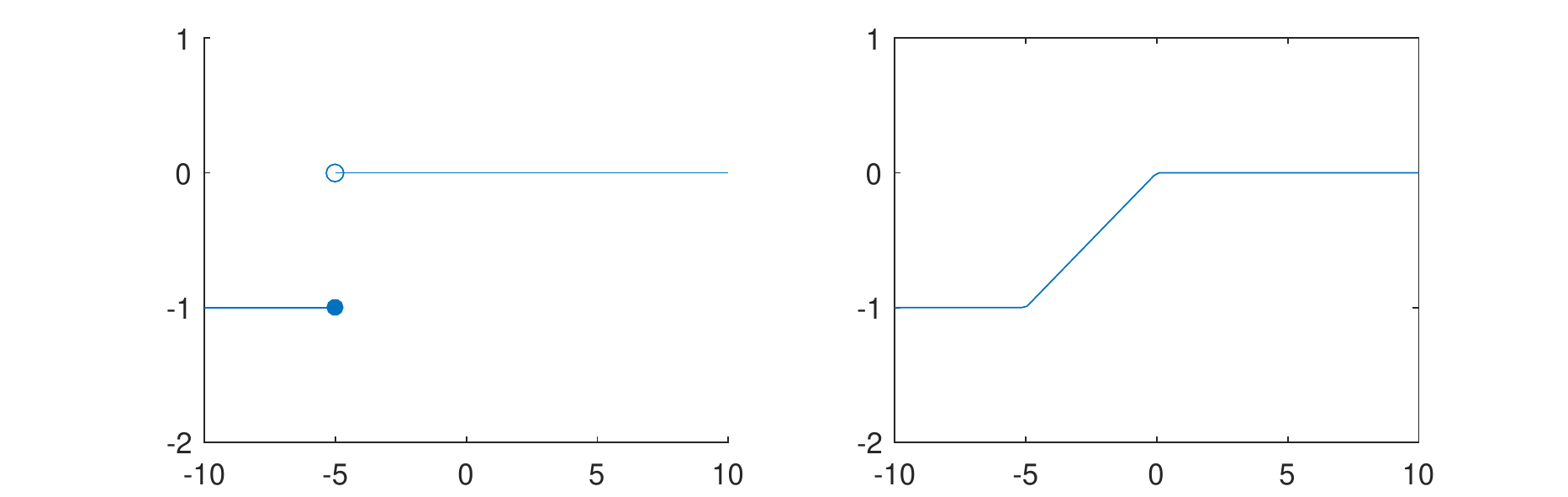}
\caption{The left plot shows the step loss function (\ref{eqn:step}), and the right plot shows the hinge loss function (\ref{eqn:hinge}).  The $x$-axis on both plots is $100\cdot(x_{T+n}/\tilde{x}_0-1)$, which is percent reduction in body weight.}
\label{fig:varphi_raw}
\end{figure}

We optimize the weight loss program using our ABMA algorithm to implement the second and third steps of our behavioral analytics framework.  This requires choosing a loss function for each individual, and Figure \ref{fig:varphi_raw} shows three choices that we considered.  These three losses make varying tradeoffs between achieving the primary health outcome of number of individuals with clinically significant weight loss (i.e., 5\% weight loss) at the end of the program versus the secondary health outcome of maximizing weight loss of individuals who were not able to achieve 5\% weight loss. The first choice of a loss function is the \emph{step loss}, which is given by
\begin{equation}
\label{eqn:step}
\varphi = \begin{cases} -1, &\text{if } x_{T+n} \leq 0.95\cdot \tilde{x}_0\\
0, &\text{otherwise}\end{cases}
\end{equation}
This discontinuous choice of a loss function gives minimal loss to 5\% or more reduction in body weight and maximal loss to less than 5\% reduction in body weight.  The second choice of a loss function is the \emph{hinge loss}, which is given by
\begin{equation}
\label{eqn:hinge}
\varphi = \begin{cases} -1, &\text{if } x_{T+n} < 0.95\cdot \tilde{x}_0\\
-0.2\cdot (x_{T+n}/\tilde{x}_0 - 1), &\text{if } 0.95\cdot \tilde{x}_0\leq x_{T+n}\leq \tilde{x}_0\\
0, &\text{if } x_{T+n} > \tilde{x}_0\end{cases}
\end{equation}
This continuous choice of a loss function gives minimal loss to 5\% or more reduction in body weight, maximal loss to less than 0\% reduction in body weight, and an intermediate loss for intermediate reductions in body weight.  The third choice of a loss function is the \emph{time-varying hinge loss}, which is given by 
\begin{equation}
\label{eqn:timvar_hinge}
\varphi = \begin{cases} -1, &\text{if } \frac{x_{T+n}}{\tilde{x}_0} <0.95  - \frac{0.05}{\log T}\\
10(\log T)(\frac{x_{T+n}}{\tilde{x}_0} - (0.95  + \frac{0.05}{\log T})), &\text{if } 0.95  - \frac{0.05}{\log T}\leq \frac{x_{T+n}}{\tilde{x}_0} \leq 0.95  + \frac{0.05}{\log T}\\
0, &\text{if } \frac{x_{T+n}}{\tilde{x}_0} > 0.95  + \frac{0.05}{\log T}\end{cases}
\end{equation}
Much like the hinge loss (\ref{eqn:hinge}), it promotes intermediate amounts of weight loss that might not meet the 5\% threshold of clinically significant weight loss. However, as more data is collected it approaches the step loss (\ref{eqn:step}) to reflect a higher degree of confidence in the estimated parameters. Thus, this choice of loss function can be considered an intermediate between the hinge (\ref{eqn:hinge}) and step (\ref{eqn:step}) losses. There is one computational note.  Since these losses are non-decreasing, we can modify Step 4 of the ABMA algorithm to instead minimize the body weight of each individual and then compose the body weight with the loss function.

%x = linspace(-10,10,100);
%subplot(121); line([-10 -5], [-1 -1]); line([-5 10], [0 0])
%axis([-10 10 -2 1])
%hold on; plot(-5,-1,'.','MarkerSize',20); plot(-5,0,'o','MarkerSize',6,'Color',[0 0.4470 0.7410]);
%subplot(122); plot(x, -(x <= -5) + (0.2*x).*(x > -5 & x <= 0))
%axis([-10 10 -2 1])

For the purpose of comparing various program designs through simulations, we considered three additional designs for the weight loss program.  We used an adaptive heuristic to design the weight loss program: Clinical visits were scheduled towards the end of the treatment at least one week apart, with more visits given to individuals who were closer to meeting the weight loss goal of a 5\% weight reduction based on their latest observed weight, and step goals were set to be a 10\% increase over a linear moving average of the individual's observed step count over the prior week. The second design was a ``do nothing'' plan where individuals were given exactly one clinical visit after two weeks, and their step goals were a constant 10,000 steps each day.  The third design was the original design of the mDPP trial: Clinical visits were scheduled on predetermined days during the treatment after 2, 4, 6, 19, 14, 18, and 20 weeks of the trials. The first two weeks of this design did not contain any clinical visits or exercise goals but instead served as an initialization period. After the first two weeks, exercise goals increased 20\% each week, starting with a 20\% increase over the average number of steps taken by individuals during the the two week initialization period. The exercise goals were capped at a maximum of 12,000 steps a day. Since the adaptive heuristic and ABMA are both adaptive, we recalculated both at the beginning of each month of the treatment and allowed both adaptive methods a 2 week initialization time similar to the mDPP trial.

\subsection{Simulation Comparison}

We compared the six different program designs using simulations of a weight loss program with a five month duration and with 30 individuals participating. Each individual in the simulation followed our behavioral model, and the parameters corresponding to the behavioral model for each individual were chosen to be those estimated by computing the MLE using the data from the mDPP trial. Since we also wanted to test how these different designs account for missing data, we assumed that the data available to each algorithm would be limited to days of the mDPP study where a particular individual reported their weight and steps.  Since the adaptive heuristic and our behavioral analytics framework are both adaptive, we recalculated the program design at the beginning of each month of the program (by re-runnning the heuristic calculations and rerunning the ABMA algorithm) and allowed both adaptive methods a two-week initialization time similar to the design of the program in the mDPP trial. All simulations were run using MATLAB on a laptop computer with a 2.4GHz processor and 16GB RAM. The Gurobi solver \citep{gurobi} in conjunction with the CVX toolbox for MATLAB \citep{cvx} were used to perform the initial estimation of the individual parameters, compute designs for the weight loss program, and perform simulations of each design.

\begin{figure}
\centering
\includegraphics[trim={0.3in 3.5in 0in 3.5in},clip,scale=0.8]{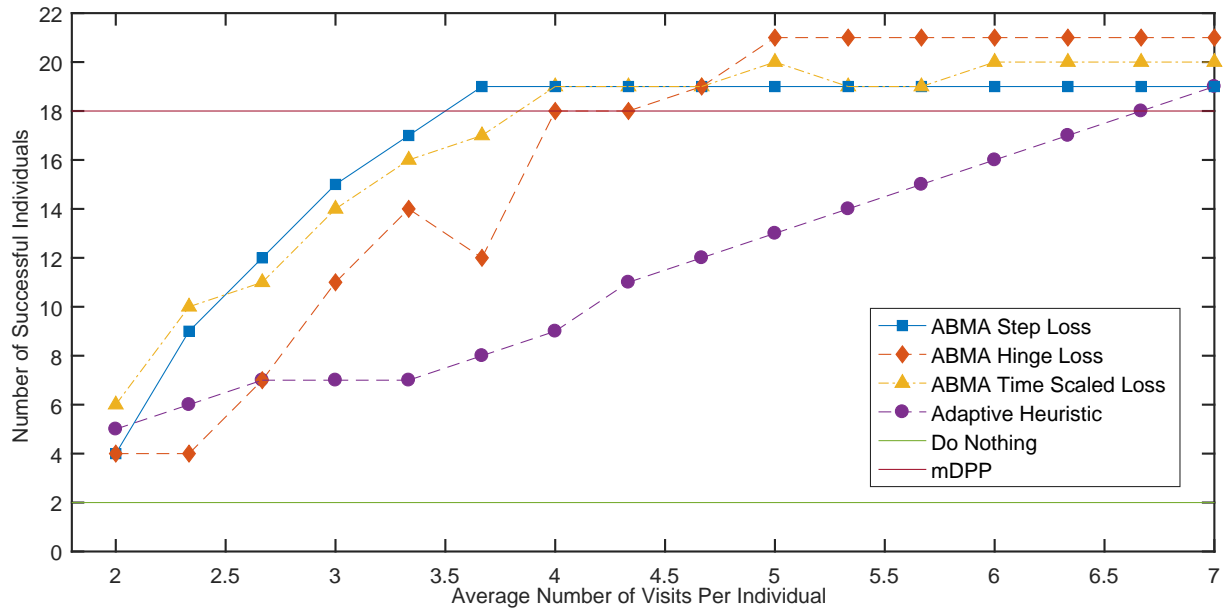}
\caption{Comparison of different program design methods with respect to number of successful individuals (i.e., lost 5\% or more body weight)}
\label{fig:num_suc_pat}
\end{figure}

Figure \ref{fig:num_suc_pat} compares the primary outcome of interest to clinicians, which is the number of individuals that achieve clinically significant weight loss (i.e., 5\% or more reduction in body weight).  We repeated the simulations for our behavioral analytics framework and the adaptive heuristic under different constraints on the total number of clinical visits that could be allocated to individuals. The $x$-axis of Figure \ref{fig:num_suc_pat} is the average number of clinical visits provided to individuals.  The horizontal line at 18 is the number of individuals who achieved 5\% weight loss in the actual mDPP trial, in which each individual received 7 clinical visits.  Figure \ref{fig:num_suc_pat} shows that all forms of behavioral analytics program and adaptive heuristic program designs outperform the ``do nothing'' policy.  Furthermore, our behavioral analytics approach and the adaptive heuristic achieve results comparable to the original mDPP program design but with significantly less resources (i.e., less clinical visits). The simulations predict that using our behavioral analytics approach in which ABMA has a step (\ref{eqn:step}) or time-varying hinge loss (\ref{eqn:timvar_hinge}) to design the weight loss program can provide health outcomes comparable to current clinical practice while using only 40-60\% of the resources (i.e., clinical visits) of current practice. In contrast, the adaptive heuristic would require 80-95\% of resources (i.e., clinical visits) to attain health outcomes comparable to current clinical practice. This suggests that appropriate choice of the loss function for our ABMA algorithm, as part of behavioral analytics approach, to personalize the design of a weight loss program could increase capacity (in terms of the number of individuals participating in the program for a fixed cost) by up to 60\%, while achieving comparable health outcomes.

\begin{figure}
\centering
\includegraphics[trim={0.3in 3.5in 0in 3.5in},clip,scale=0.8]{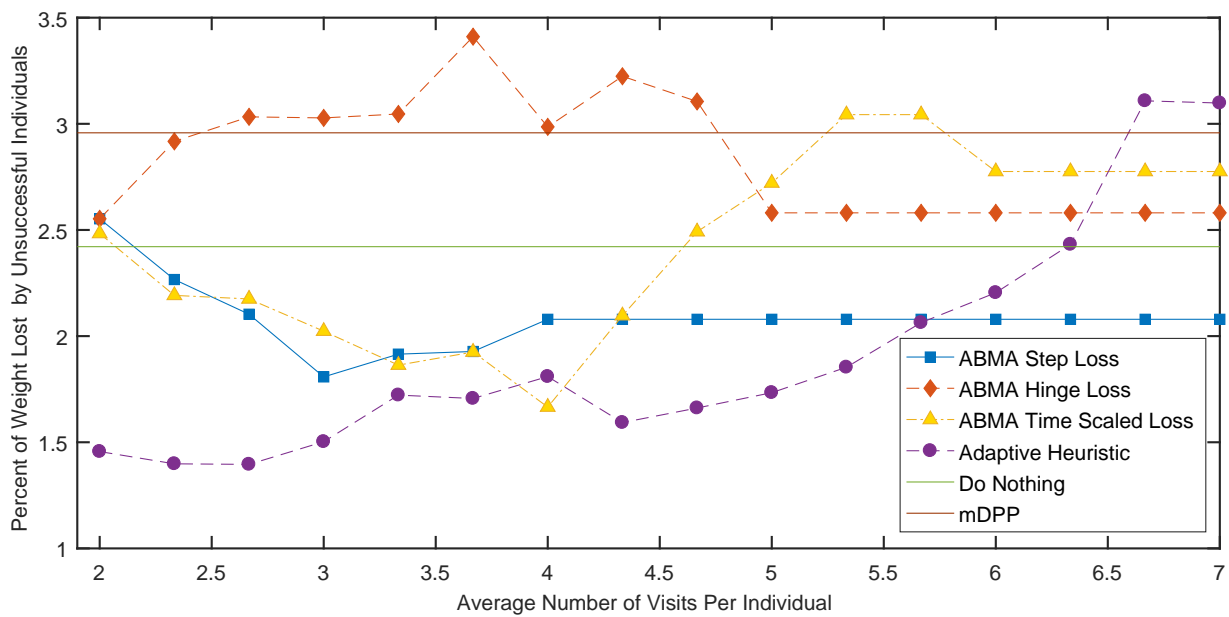}
		\caption{Comparison of different program design methods with respect to the percent weight lost by unsuccessful individuals (i.e., lost less than 5\% body weight)}
		\label{fig:cvar}
\end{figure}

Figure \ref{fig:cvar} compares the different program designs using a secondary outcome of interest to clinicians of the average amount of weight loss of individuals who did not successfully achieve 5\%  weight loss.  The original treatment plan of mDPP and the ``do nothing'' treatment plan slightly outperform the adaptive program designs at certain clinical visit budgets. This effect however is mainly due to these static plans not identifying individuals who are on the cusp of achieving 5\% weight loss but might still achieve around 3-4\% weight loss, while both adaptive program designs allocate clinical visits to these individuals and ensure they reach the weight loss goal of 5\% weight loss. Restated, the lower weight loss of unsuccessful individuals under the behavioral analytics treatment plans is an artifact of the improved success rate of the behavioral analytics plans in helping individuals achieve 5\% weight loss. This effect is further exemplified in the region of between an average of 2.8-4.2 visits per individual, where we see that individuals who were not successful in achieving 5\% weight loss in the behavioral analytics treatment plans on average lose more weight then those under the heuristic policy, while in the region of an average of 5.6-7 visits per individual we see that individuals under the heuristic treatment lose more weight. The effect in this last region is mainly due to the behavioral model used in our behavioral analytics framework, which is more effective at identifying the individuals who would most benefit from additional clinical visits. Therefore, more resources are spent on individuals who could potentially reach their 5\% weight loss goal, while the adaptive heuristic uses these resources in a less effective manner. 

Figures \ref{fig:num_suc_pat} and \ref{fig:cvar} demonstrate a tradeoff between the primary and secondary outcomes, and the various loss functions provide different tradeoffs. Note the line for the step loss (\ref{eqn:step}) is the first to achieve a primary outcome comparable to that of the original mDPP trial while fluctuating relatively little in terms of the secondary outcome. This matches intuition that the step loss function (\ref{eqn:step}) is focused on ensuring individuals achieve 5\% weight loss while not being concerned with their final weight. On the other hand, the line for the hinge loss (\ref{eqn:hinge}) lags behind the other behavioral analytics policies in achieving comparable primary outcomes to the mDPP trial while having an extremely effective secondary outcome. These results follow our intuition that this loss favors intermediate weight loss over achieving clinically significant weight loss. Finally, the line corresponding to the time-varying hinge loss function (\ref{eqn:timvar_hinge}) has a clear transition at an average of 5 visits per individual from favoring the primary outcome to the secondary outcome. This behavior indicates that using such time scaling leads to interventions that focus on primary outcomes when resources are constrained but also accounts for secondary outcomes when resources are less scarce. Such behavior may be useful for implementing a behavioral analytics approach when the relative abundance of resources is not known \emph{a priori}.  

\subsection{Computational Performance and Sensitivity Analysis}
\begin{figure}
\centering
\includegraphics[trim={0.3in 3.5in 0in 3.5in},clip,scale=0.8]{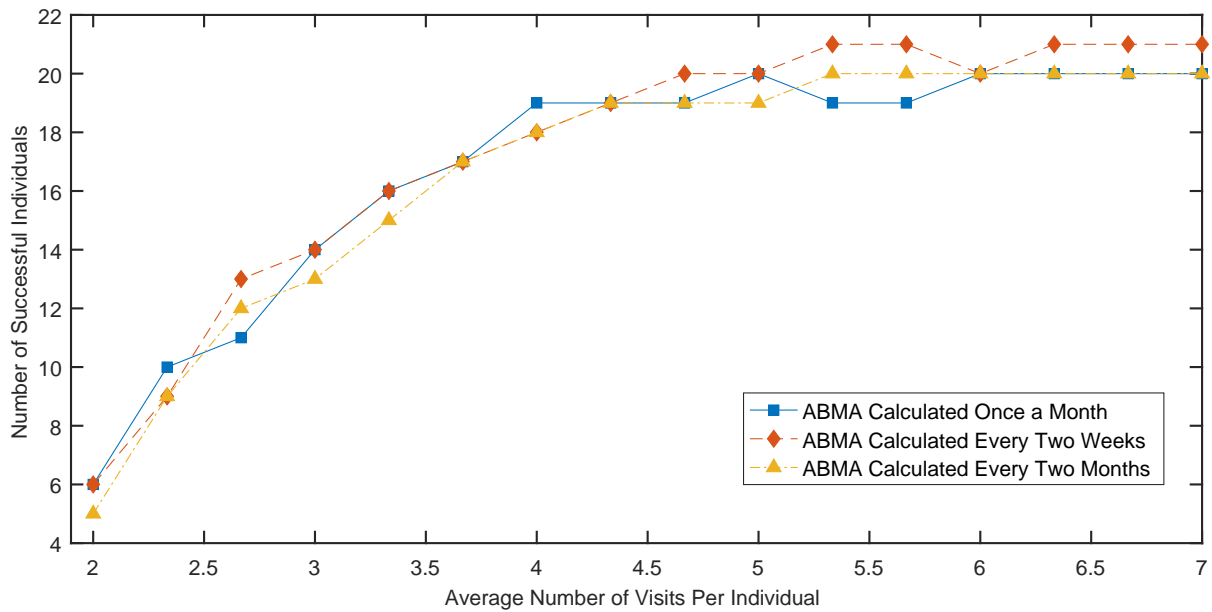}
\caption{ Comparison of calculation schedules and their effects on the number of successful individuals (i.e., lost 5\% or more body weight)}
\label{fig:sched_numpat}
\end{figure}

\begin{figure}
\centering
\includegraphics[trim={0.3in 3.5in 0in 3.5in},clip,scale=0.8]{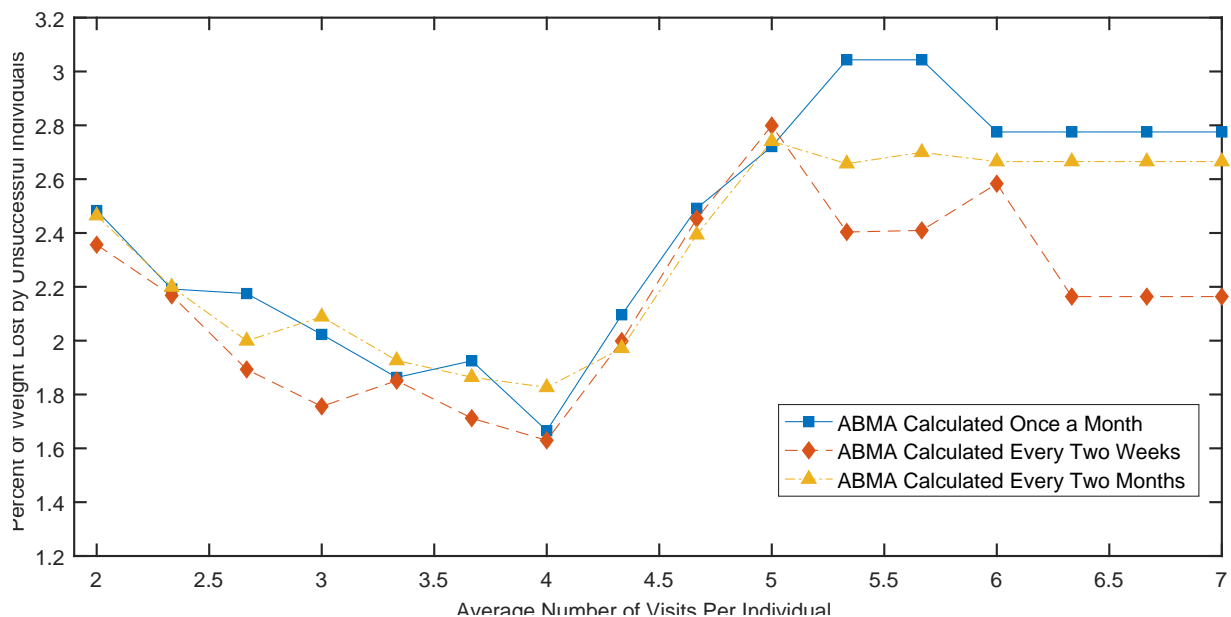}
\caption{Comparison of calculation schedules and their effects on the weight lost by unsuccessful individuals (i.e., lost less than 5\% body weight)}\label{fig:sched_cvar}
\end{figure}

The simulation experiments assumed that treatment plans were updated at the beginning of each month by re-running the second and third steps of our behavioral analytics approach (through applying the ABMA algorithm), and so we conducted a sensitivity analysis to examine the effect of updating the program design more or less frequently. Figures \ref{fig:sched_numpat} and \ref{fig:sched_cvar} compare the health outcomes of using a program designed by our behavioral analytics framework with a time-varying hinge loss (\ref{eqn:timvar_hinge}), where the treatment was recalculated once every two weeks, once a month, and once every two months. These results show that recomputing the treatment plan with lesser or higher frequency does not significantly impact the efficacy of the resulting treatment. This indicates that for practical implementation, the statistical convergence rate of estimated parameters in our behavioral model is sufficiently fast that it would suffice to rerun the weight loss program design algorithm at most once a month.

We also conducted time-benchmarks for the sub-problems involved in computing 2SSA, MAP, and ABMA, which are the algorithms comprising the second and third steps of our implementation o a behavioral analytics framework.  The results of the time-benchmarks are summarized in Tables \ref{tabel:int_array},\ref{tabel:pdf_calc}, and \ref{tabel:knap_calc} in the appendix. On average, solving all sub-problems took 17s per individual. This is promising for practical implementation, particularly because each sub-problem calculation can be performed in parallel for each size constraint of the clinical visit schedule (from 1 to 7 visits).  The results show that computation time increases with respect to the number of available visits and data available in the treatment plan calculation. However, the calculation times still remain below 30s on average per individual for each step of the program calculation. This would imply that our methodology for weight loss program calculation is suitable for large scale program design since the program design would be updated at most once every month.

\section{Conclusion}
In this paper, we develop a {\em behavioral analytics} framework for multi-agent systems in which a single coordinator provides behavioral or financial incentives to a large number of myopic agents. Our framework is applicable in a variety of settings of interest to the operations research community, including the design of demand-response programs for electricity consumes,  the personalized design of a weight loss program, and adaptive logistics allocation for franchises.   The framework we develop involves the definition of a behavioral model, the estimation of model parameters, and the optimization of incentives.   We show (among other results) that under mild assumptions, the incentives computed by our approach converge to the optimal incentives that would be computed knowing full information about the agents.  We evaluated our approach for personalizing the design of a weight loss program, and showed via simulation that our approach can improve outcomes with reduced treatment cost.    
 \begin{APPENDICES}
 \section{Complete MILP Formulation of MLE Problem}
\begin{align}
\min\ &  \frac{\sqrt{2}}{\sigma_1}\sum\limits_{i=1}^{n_w} \zeta_{w,t_i} + \frac{\sqrt{2}}{\sigma_2}\sum\limits_{i=1}^{n_u} \zeta_{u,t_i} \\
\text{s.t. }&-\zeta_{w,t_i} \leq \tilde{w}_{t_i} - w_{t_i} \leq \zeta_{w,t_i} &\forall 1\leq t_i \leq n_w\\
& -\zeta_{u,t_i} \leq \tilde{u}_{\tau_i} - u_{\tau_i} \leq \zeta_{u,t_i} &\forall1\leq t_i \leq n_u\\
&-b(aw_1 +bu_0 + f_0 + k) -r(u_0-u_b) = 0 \\
&-(aw_1 + bu_0 + f_0 + k) - r(f_0-F_{b,1}) = 0\\
&-b(aw_t +bu_t + f_t + k) -r(u_t-u_b) = 0 & \forall 1 \leq t \leq m\\
&-(aw_t + bu_t + f_t + k) - r(f_t-f_{b,t})= 0 & \forall 1 \leq t \leq m\\
&-2b(aw_t + bu_t + f_t +k) -2r(u_t-u_b) + \lambda_{1,t} = 0  &\forall m\leq t \leq n\\
&-2(aw_t +bu_t +f_t +k) -2(f_t-f_{b,t}) =0 &\forall m\leq t \leq n\\
&(g_t - \epsilon) - Mx_{1,t} \leq u_t \leq g_t - \epsilon + M(1 - x_{1,t}) &\forall m\leq t \leq n\\
&(g_t - \epsilon) - M(1- x_{2,t}) \leq u_t \leq g_t + \epsilon + M(1- x_{2,t}) &\forall m\leq t \leq n\\
&(g_t + \epsilon) - M(1 - x_{3,t}) \leq u_t \leq g_t + \epsilon + Mx_{3,t} &\forall m\leq t \leq n\\
&p_t -M(1-x_{t,1}) \leq \lambda_{1,t} \leq M(1-x_{3,t})  &\forall m\leq t \leq n\\
&0 \leq \lambda_{1,t} \leq p_t &\forall m\leq t \leq n\\
&x_{t,1} + x_{t,2} + x_{t,3} = 1 &\forall m\leq t \leq n\\
&p_{t+1} \geq \gamma p_t + \delta d_{t+1}  &\forall m\leq t \leq n\\
&p_{t+1} \leq \gamma p_t + \delta d_{t+1} + M(1 - x_{1,t}) &\forall m\leq t \leq n\\
&p_{t+1} \geq \gamma p_t + \delta d_{t+1} + \mu - M x_{1,t} &\forall m\leq t \leq n\\
&p_{t+1} \leq \gamma p_t + \delta d_{t+1} + \mu &\forall m\leq t \leq n\\
&F_{b,t+1} = (1-\alpha)F_{b,t} + \alpha f_{b,t} & \forall t \\
&f_{b,t+1} = \gamma (f_{b,t} - F_{b,t}) + F_{b,t} - \beta d_t & \forall t \\
&x_{t+1,1} \geq x_{t,1} - d_{t+1} - \mathds{1}(g_{t+1} - g_t < 0) &\forall m\leq t \leq n\\
&x_{t+1,2} \leq x_{t,2} - d_{t+1} + \mathds{1}(g_{t+1} - g_t < 0)&\forall m\leq t \leq n\\
&x_{t+1,3} \leq x_{t,3} - d_{t+1} + \mathds{1}(g_{t+1} - g_t < 0)&\forall m\leq t \leq n\\
&w_t,u_t,f_t,F_{b,t},f_{b,t},p_{t} \geq 0 & \forall t \\
&x_{t,1},x_{t,2},x_{t,3} \in \mathbb{B} & \forall t
\end{align}
 \section{Complete MILP Formulation of MAP Problem}
\begin{align}
\min\ &  \frac{\sqrt{2}}{\sigma_1}\sum\limits_{i=1}^{n_w} \zeta_{w,t_i} + \frac{\sqrt{2}}{\sigma_2}\sum\limits_{i=1}^{n_u} \zeta_{u,t_i} - \sum\limits_{x \in \theta} \sum\limits_{i=1}^{m_x}z_{i,x} \log \phi_i^x \\
\text{s.t. }&-\zeta_{w,t_i} \leq \tilde{w}_{t_i} - w_{t_i} \leq \zeta_{w,t_i} &\forall 1\leq t_i \leq n_w\\
& -\zeta_{u,t_i} \leq \tilde{u}_{\tau_i} - u_{\tau_i} \leq \zeta_{u,t_i} &\forall1\leq t_i \leq n_u\\
&-b(aw_1 +bu_0 + f_0 + k) -r(u_0-u_b) = 0 \\
&-(aw_1 + bu_0 + f_0 + k) - r(f_0-F_{b,1}) = 0\\
&-b(aw_t +bu_t + f_t + k) -r(u_t-u_b) = 0 & \forall 1 \leq t \leq m\\
&-(aw_t + bu_t + f_t + k) - r(f_t-f_{b,t})= 0 & \forall 1 \leq t \leq m\\
&-2b(aw_t + bu_t + f_t +k) -2r(u_t-u_b) + \lambda_{1,t} = 0  &\forall m\leq t \leq n\\
&-2(aw_t +bu_t +f_t +k) -2(f_t-f_{b,t}) =0 &\forall m\leq t \leq n\\
&(g_t - \epsilon) - Mx_{1,t} \leq u_t \leq g_t - \epsilon + M(1 - x_{1,t}) &\forall m\leq t \leq n\\
&(g_t - \epsilon) - M(1- x_{2,t}) \leq u_t \leq g_t + \epsilon + M(1- x_{2,t}) &\forall m\leq t \leq n\\
&(g_t + \epsilon) - M(1 - x_{3,t}) \leq u_t \leq g_t + \epsilon + Mx_{3,t} &\forall m\leq t \leq n\\
&p_t -M(1-x_{t,1}) \leq \lambda_{1,t} \leq M(1-x_{3,t})  &\forall m\leq t \leq n\\
&0 \leq \lambda_{1,t} \leq p_t &\forall m\leq t \leq n\\
&x_{t,1} + x_{t,2} + x_{t,3} = 1 &\forall m\leq t \leq n\\
&p_{t+1} \geq \gamma p_t + \delta d_{t+1}  &\forall m\leq t \leq n\\
&p_{t+1} \leq \gamma p_t + \delta d_{t+1} + M(1 - x_{1,t}) &\forall m\leq t \leq n\\
&p_{t+1} \geq \gamma p_t + \delta d_{t+1} + \mu - M x_{1,t} &\forall m\leq t \leq n\\
&p_{t+1} \leq \gamma p_t + \delta d_{t+1} + \mu &\forall m\leq t \leq n\\
&F_{b,t+1} = (1-\alpha)F_{b,t} + \alpha f_{b,t} & \forall t \\
&f_{b,t+1} = \gamma (f_{b,t} - F_{b,t}) + F_{b,t} - \beta d_t & \forall t \\
&x_{t+1,1} \geq x_{t,1} - d_{t+1} - \mathds{1}(g_{t+1} - g_t < 0) &\forall m\leq t \leq n\\
&x_{t+1,2} \leq x_{t,2} - d_{t+1} + \mathds{1}(g_{t+1} - g_t < 0)&\forall m\leq t \leq n\\
&x_{t+1,3} \leq x_{t,3} - d_{t+1} + \mathds{1}(g_{t+1} - g_t < 0)&\forall m\leq t \leq n\\
&z_{i,x}h_{lb,i}^x \leq x_i \leq z_{i,x}h_{ub,i}^x &\forall x \forall i\\
&\sum\limits_{i=1}^{m_x} z_{i,x} = 1 &\forall x \forall i\\
&\sum\limits_{i=1}^{m_x} x_i = x &\forall x\forall i\\
&z_{i,x} \in \mathbb{B} &\forall x \forall i \\
&w_t,u_t,f_t,F_{b,t},f_{b,t},p_{t} \geq 0 & \forall t \\
&x_{t,1},x_{t,2},x_{t,3} \in \mathbb{B} & \forall t
\end{align}

\section{Complete MILP Formulation of Personalized Treatment Plan Design Problem}
\begin{align}
\min\ & \; w_n \\
\text{s.t. }&F_{b,t+1} = (1-\alpha)F_{b,t} + \alpha f_{b,t} &\forall t>T\\
&f_{b,t+1} = \gamma (f_{b,t} - F_{b,t}) + F_{b,t} - \hat{\beta} d_t &\forall t>T\\
&p_{t+1} = \gamma p_t + \hat{\delta} d_t + \hat{\mu} (1- x_{1,t}) &\forall t>T\\ 
&d_t \leq  \mathds{1}\left[\text{mod}(t,7) = 1\right] & \forall t > T\\
&d_t \leq 1- d_\tau \quad \forall \tau> T,  \tau +1 \leq t \leq \tau + 6 & \forall t > T \\
&-2b(aw_t + bu_t + f_t +k) -2r(u_t-u_b) + \lambda_{1,t} + \lambda_{4,t} = 0 & \forall t > T\\
&-2(aw_t +bu_t +f_t +k) -2(f_t-f_{b,t}) + \lambda_{3,t}=0 & \forall t > T\\
&(g_t - \epsilon) - Mx_{1,t} \leq u_t \leq g_t - \epsilon + M(1 - x_{1,t})& \forall t > T \\
&(g_t - \epsilon) - M(1- x_{2,t}) \leq u_t \leq g_t + \epsilon + M(1- x_{2,t}) & \forall t > T\\
&(g_t + \epsilon) - M(1 - x_{3,t}) \leq u_t \leq g_t + \epsilon + Mx_{3,t}& \forall t > T\\
&p_t -M(1-x_{t,1}) \leq \lambda_{1,t} \leq M(1-x_{3,t}) & \forall t > T\\
&0 \leq f_t \leq M(1-x_{f,t}); 0 \leq u_t \leq M(1-x_{u,t})& \forall t > T\\
&0 \leq \lambda_{3,t} \leq M x_{f,t} & \forall t > T\\
&0 \leq \lambda_{4,t} \leq M x_{u,t}& \forall t > T\\
&0 \leq \lambda_{1,t} \leq p_t & \forall t > T\\
&x_{t,1} + x_{t,2} + x_{t,3} = 1 &\forall t > T\\
&g_{t+1} - g_t \leq M(1 - g_{ind,t}) &\forall t > T\\
&g_{t+1} - g_t \geq -M g_{ind,t} &\forall t > T\\
&x_{t+1,1} \geq x_{t,1} - d_{t+1} - g_{ind,t}&\forall t > T\\
&x_{t+1,2} \leq x_{t,2} - d_{t+1} + g_{ind,t}&\forall t > T\\
&x_{t+1,3} \leq x_{t,3} - d_{t+1} + g_{ind,t}&\forall t > T \\
&g_{ind,t} \in \mathbb{B} \\
&x_{t,1},x_{t,2},x_{t,3}, x_{f,t}, x_{u,t} \in \mathbb{B} \\
&d_t = \bar{d}_t; g_t = \bar{g}_t; w_t = \hat{w}_t; u_t = \hat{u}_t; f_t = \hat{f}_t; \theta_t = \hat{\theta}_t  & \forall t
\end{align}

\section{Benchmarking Performance Tables}

 \begin{table}[!htb]
\centering
\tiny
\begin{tabular}{|c|c|c c c c c |c|}
  \hline
  \multicolumn{8}{|c|}{Average Runtimes for Candidate Treatment Plan Calculation (in seconds)} \\
  \hline
	&\multicolumn{7}{c|}{Date of Calculation During the Program} \\
  &&15&30&60&90&120&Average \\
	\hline
	\multirow{18}{*}{Visit Budget} &
54&21.715&9.3363&8.8984&9.5618&10.003&11.903 \\
&63&21.715&7.7624&14.714&18.736&16.692&15.924 \\
&72&21.715&8.8894&21.239&26.186&16.091&18.824 \\
&81&21.715&11.408&43.147&28.112&23.284&25.533 \\
&90&21.715&10.274&23.935&13.421&33.6&20.589 \\
&99&21.715&9.3366&26.251&14.16&19.855&18.263 \\
&108&21.715&9.4194&24.208&12.774&21.757&17.975 \\
&117&21.715&9.9031&22.962&9.9226&31.865&19.273 \\
&126&21.715&10.217&21.749&10.199&28.664&18.509 \\
&135&21.715&9.973&16.307&23.129&18.893&18.003 \\
&144&21.715&11.208&14.626&14.851&8.9203&14.264 \\
&153&21.715&12.978&13.913&13.892&8.0501&14.11 \\
&162&21.715&12.403&14.674&14.434&10.965&14.838 \\
&171&21.715&13.305&12.102&11.585&23.78&16.497 \\
&180&21.715&14.879&12.116&11.964&18.473&15.829 \\
&189&21.715&11.731&11.835&12.221&21.715&15.843 \\
\hline
&Average &21.715&10.814&18.917&15.322&19.538&17.261 \\
  \hline
\end{tabular}
\caption{}
\label{tabel:int_array}
\end{table}

\begin{table}[!htb]
\centering
\tiny
\begin{tabular}{|c|c|c c c c c |c|}
  \hline
  \multicolumn{8}{|c|}{Average Runtimes for MAP Calculation (in seconds)} \\
  \hline
	&\multicolumn{7}{c|}{Date of Calculation During the Program} \\
  &&15&30&60&90&120&Average \\
	\hline
	\multirow{18}{*}{Visit Budget} &
54&16.365&11.416&9.537&11.479&7.6851&11.296 \\
&63&16.365&13.249&12.816&17.157&11.493&14.216 \\
&72&16.365&11.913&14.15&13.846&12.884&13.832 \\
&81&16.365&17.448&16.11&11.936&17.237&15.819 \\
&90&16.365&12.251&9.9218&8.9877&16.463&12.798 \\
&99&16.365&9.9473&10.59&9.3482&13.666&11.983 \\
&108&16.365&10.412&10.199&9.377&14.746&12.22 \\
&117&16.365&10.696&9.5879&8.9027&23.187&13.748 \\
&126&16.365&10.36&10.737&9.6123&23.524&14.12 \\
&135&16.365&9.7085&10.241&12.753&18.221&13.458 \\
&144&16.365&9.3022&13.258&12.135&11.654&12.543 \\
&153&16.365&8.5183&11.297&11.527&11.152&11.772 \\
&162&16.365&9.4638&8.4174&11.279&14.812&12.068 \\
&171&16.365&8.0878&11.092&9.8865&18.003&12.687 \\
&180&16.365&8.8929&7.5454&6.5939&12.739&10.427 \\
&189&16.365&7.679&7.6401&6.9321&16.365&10.996 \\
\hline
&Average &16.365&10.584&10.821&10.735&15.239&12.749 \\
  \hline
\end{tabular}
\caption{}
\label{tabel:pdf_calc}
\end{table}

\begin{table}[!htb]
\centering
\tiny
\begin{tabular}{|c|c|c c c c c |c|}
  \hline
  \multicolumn{8}{|c|}{Average Runtimes for Knapsack Calculation (in seconds)} \\
  \hline
	&\multicolumn{7}{c|}{Date of Calculation During the Program} \\
  &&15&30&60&90&120&Average \\
	\hline
	\multirow{18}{*}{Visit Budget} &
54&0.21791&0.19278&0.14949&0.1666&0.19679&0.18471 \\
&63&0.19295&0.16084&0.19613&0.23577&0.22246&0.20163 \\
&72&0.21302&0.2516&0.17474&0.18052&0.3125&0.22648 \\
&81&0.25345&0.3542&0.18128&0.17449&0.25997&0.24468 \\
&90&0.21698&0.17154&0.168&0.16942&0.21394&0.18798 \\
&99&0.17072&0.17607&0.15058&0.15699&0.20012&0.1709 \\
&108&0.17086&0.17178&0.17235&0.15172&0.19029&0.1714 \\
&117&0.17958&0.17073&0.16465&0.16171&0.20232&0.1758 \\
&126&0.18599&0.24335&0.17326&0.16726&0.20706&0.19538 \\
&135&0.20403&0.16143&0.16626&0.16681&0.19607&0.17892 \\
&144&0.17733&0.14174&0.16112&0.1618&0.18708&0.16581 \\
&153&0.18841&0.15411&0.17753&0.17965&0.18822&0.17758 \\
&162&0.17892&0.2167&0.18539&0.17419&0.21756&0.19455 \\
&171&0.1771&0.29787&0.21357&0.2248&0.21789&0.22625 \\
&180&0.20535&0.18508&0.22285&0.16565&0.19392&0.19457 \\
&189&0.18554&0.2104&0.19875&0.17776&0.20053&0.1946 \\
\hline
&Average &0.19488&0.20376&0.1785&0.17595&0.21292&0.1932 \\
  \hline
\end{tabular}
\caption{}
\label{tabel:knap_calc}
\end{table}

\end{APPENDICES}
\newpage
\ACKNOWLEDGMENT{The authors gratefully acknowledge the support of NSF Award CMMI-1450963, UCSF Diabetes Family Fund for Innovative Patient Care-Education and Scientific Discovery Award, K23 Award (NR011454), and the UCSF Clinical and Translational Science Institute (CTSI) as part of the Clinical and Translational Science Award program funded by NIH UL1 TR000004.}

\bibliographystyle{ormsv080} % outcomment this and next line in Case 1
\bibliography{soba} % if more than one, comma separated

%%%%%%%%%%%%%%%%%%
\end{document}